\documentclass[onefignum,onetabnum]{siamart220329}

\usepackage{lipsum}
\usepackage{subcaption}
\usepackage{amsfonts}
\usepackage{graphicx}
\usepackage{epstopdf}
\usepackage{appendix}
\usepackage{url}
\usepackage{algpseudocode}
\ifpdf
  \DeclareGraphicsExtensions{.eps,.pdf,.png,.jpg}
\else
  \DeclareGraphicsExtensions{.eps}
\fi

\usepackage{enumitem}
\setlist[enumerate]{leftmargin=.5in}
\setlist[itemize]{leftmargin=.5in}

\newsiamremark{remark}{Remark}
\newsiamremark{hypothesis}{Hypothesis}
\crefname{hypothesis}{Hypothesis}{Hypotheses}
\newsiamthm{claim}{Claim}

\headers{Faster Legendre-Chebyshev transforms}{M. Mortensen}

\title{An Example Article\thanks{Submitted to the editors DATE.
\funding{This work was funded by the Fog Research Institute under contract no.~FRI-454.}}}

\author{Mikael Mortensen\thanks{University of Oslo 
  (\email{mikaem@math.uio.no}, \url{https://www.mn.uio.no/math/english/people/aca/mikaem/}).}
}
\usepackage{amsopn}
\DeclareMathOperator{\diag}{diag}

\usepackage{enumitem}
\usepackage{array}
\usepackage{diagbox}
\usepackage{xfrac}
\usepackage{listings}
\usepackage{subcaption}

\usepackage{tikzit}

\tikzstyle{new style 0}=[fill=white, draw=black, shape=rectangle, minimum width=10cm, minimum height=4cm, tikzit fill=white, tikzit draw=black, tikzit shape=rectangle]
\tikzstyle{dot}=[fill=black, draw=black, shape=circle, tikzit fill=black, tikzit draw=black, inner sep=1pt]

\tikzstyle{new edge style 0}=[-, fill=blue, tikzit fill=blue, draw=blue, tikzit draw=blue]
\tikzstyle{new edge style 1}=[{|-|}]
\tikzstyle{arrow}=[->]
\tikzstyle{thick}=[line width=1.5pt, ->]
\tikzstyle{thin}=[-, line width=0.4pt, dotted]

\ifpdf
\hypersetup{
  pdftitle={A faster multipole Legendre-Chebyshev transform},
  pdfauthor={M. Mortensen}
}
\fi

\title{A faster multipole Legendre-Chebyshev transform}

\author{Mikael Mortensen\thanks{Department of Mathematics, University of Oslo, Norway
  (\email{mikaem@math.uio.no}).}}

\date{\today}

\algnewcommand \Input[3]  {\item\hspace*{0ex}\textbf{input}  #1 : #2 \Comment{#3}}
\algnewcommand \Output[3] {\item\hspace*{0ex}\textbf{output} #1 : #2 \Comment{#3}}
\algnewcommand \WorkArray[2]  {\item\hspace*{0ex}\textbf{work array}  #1 : #2}

\begin{document}
\maketitle

\begin{abstract}
This paper describes a fast algorithm for transforming Legendre coefficients into Chebyshev coefficients, and vice versa. The algorithm is based on the fast multipole method and is similar to the approach described by Alpert and Rokhlin [SIAM J. Sci. Comput., 12 (1991)]. The main difference is that we utilise a modal Galerkin approach with Chebyshev basis functions instead of a nodal approach with a Lagrange basis. Part of the algorithm is a novel method that facilitates faster spreading of intermediate results through neighbouring levels of hierarchical matrices. This enhancement leads to a method that is approximately $20~\%$ faster to execute, due to less floating point operations. We also describe an efficient initialization algorithm that for the Lagrange basis is roughly 5 times faster than the original method for large input arrays. The described method has both a planning and execution stage that asymptotically require $\mathcal{O}(N)$ flops. The algorithm is simple enough that it can be implemented in 100 lines of vectorized Python code. Moreover, its efficiency is such that a single-threaded C implementation can transform $10^6$ coefficients in approximately 20 milliseconds on a new MacBook Pro M3, representing about 3 times the execution time of a well-planned (single-threaded) type 2 discrete cosine transform from FFTW (www.fftw.org). Planning for the $10^6$ coefficients requires approximately 50 milliseconds.
\end{abstract}

\section{Introduction}

Let $\mathbb{P}_N$ be the set of polynomials of degree less than or equal to $N$. Any function $f(x) \in \mathbb{P}_{N-1} $ defined on $ x \in [-1, 1]$ can then be expanded as
\begin{equation}
    f(x) = \sum_{j=0}^{N-1} \hat{f}^{leg}_j L_{j}(x) = \sum_{j=0}^{N-1} \hat{f}^{cheb}_j T_{j}(x), \label{eq:fx}
\end{equation}
where $\{L_j\}_{j=0}^{N-1}$ and $\{T_j\}_{j=0}^{N-1}$ are the first $N$ Legendre polynomials and Chebyshev polynomials of the first kind, respectively, and $\boldsymbol{f}^{leg} = \{\hat{f}^{leg}_j\}_{j=0}^{N-1}$ and $ \boldsymbol{f}^{cheb} = \{\hat{f}^{cheb}_j\}_{j=0}^{N-1}$ are expansion coefficients. It is the purpose of this paper to describe a fast method for computing $\boldsymbol{f}^{cheb}$ from $\boldsymbol{f}^{leg}$ and vice versa.

A direct method can be obtained by using the orthogonality of Chebyshev polynomials in the $L^2_{\omega}(-1,1)$ space, defined with the weighted inner product
\begin{equation}
    (f, v)_{\omega} = \int_{-1}^1 f v \omega dx,
\end{equation}
where the weight $\omega(x)=(1-x^2)^{-\sfrac{1}{2}}$ and $f(x)$ and $v(x)$ are real functions. We multiply the two equal sums in \eqref{eq:fx} by the test function $T_i$ and the weight $\omega$ and integrate over the domain to obtain
\begin{equation}
    \sum_{j=0}^{N-1}\left(T_i, L_{j}\right)_{\omega} \hat{f}^{leg}_j = \sum_{j=0}^{N-1} \left( T_i, T_{j}\right)_{\omega} \hat{f}^{cheb}_j, \quad i = 0, 1, \ldots, N-1.
\end{equation}
Here $(T_i, T_j)_{\omega} = c_i \pi / 2 \delta_{ij}$, where $c_0=2$, $c_i=1$ for $i>0$, $\delta_{ij}$ is the Kronecker delta-function and the diagonal matrix is denoted as $\boldsymbol{C} = \diag(\{c_i \pi/2\}_{i=0}^{N-1}) \in \mathbb{R}^{N \times N}$. The matrix on the left hand side is denoted as $a_{ij} = \left(T_i, L_{j}\right)_{\omega}$, with $ \boldsymbol{A} = (a_{ij})_{i,j=0}^{N-1} \in \mathbb{R}^{N \times N}$, such that
\begin{equation}
    \boldsymbol{A} \boldsymbol{f}^{leg} = \boldsymbol{C} \boldsymbol{f}^{cheb}, \label{eq:AfCf}
\end{equation}
and thus 
\begin{equation}
    \boldsymbol{f}^{cheb} = \boldsymbol{C}^{-1}\boldsymbol{A} \boldsymbol{f}^{leg} \quad \text{ and } \quad
    \boldsymbol{f}^{leg} = \boldsymbol{A}^{-1}\boldsymbol{C} \boldsymbol{f}^{cheb}. \label{eq:directmv}
\end{equation}
The diagonal matrix $\boldsymbol{C}$ is obviously trivial to apply to any vector, and thus we will focus on the matrix $\boldsymbol{A} $, where all nonzero elements can be represented as
\begin{equation}
    a_{ij} = \Lambda\left(\frac{j-i}{2}\right) \Lambda\left(\frac{j+i}{2}\right), \label{eq:aij}
\end{equation}
where
\begin{equation}
    \Lambda(x) = \frac{\Gamma(x+\tfrac{1}{2})}{\Gamma(x+1)}, \quad x \in \mathbb{R}^+. \label{eq:Lambda}
\end{equation}
The upper triangular matrix $\boldsymbol{A}$ can be even more easily represented diagonalwise, where the only nonzero items are
\begin{equation}
    a_{i, i+2k} =
     \Lambda(k) \Lambda(i+k), \quad \forall \, k=0, 1, \ldots, N/2-1 \text{ and }i+k<N.
     \label{eq:hik}
\end{equation}
Hence, if we define a vector $ \boldsymbol{\lambda} = \{\lambda_k\}_{k=0}^{N-1} \in \mathbb{R}^{N}$, where $ \lambda_k = \Lambda(k)$, 
then the entire matrix $\boldsymbol{A}$ can be represented using merely a single vector $\boldsymbol{\lambda}$ of storage. Also note that since every other upper diagonal is zero, the matrix can be efficiently split into odd and even parts. A direct algorithm using $\boldsymbol{A}$ in Eq. \eqref{eq:directmv} is briefly described in Sec.~\ref{sec:direct}.

The transforms between $\boldsymbol{f}^{leg}$ and $\boldsymbol{f}^{cheb}$ are useful because, although Chebyshev and Legendre expansions share similar approximation characteristics, each offers distinct advantages and drawbacks.
A Chebyshev series is orthogonal in the weighted $L^2_{\omega}(-1,1)$ space, and it is often favoured because algorithms can be accelerated through the use of cosines and fast discrete cosine transforms. The Legendre polynomials, on the other hand, are often favoured for their orthogonality in the more convenient $L^2(-1,1)$ space. However, the Legendre polynomials lack really fast transforms, which makes them slow to evaluate on a large physical mesh. The Chebyshev-Legendre Galerkin method of Shen \cite{shen96} takes advantage of the best of both worlds, making use of a transform between $\boldsymbol{f}^{leg}$ and $\boldsymbol{f}^{cheb}$, combined with fast cosine transforms.  

There are several methods available for transforming from $\boldsymbol{f}^{leg}$ to $\boldsymbol{f}^{cheb}$ and vice versa. To date the state-of-the-art in widely recognised to be the fast multipole method described by Alpert and Rokhlin \cite{alpert91}, which completes the transform in $\mathcal{O}(N)$ operations after an initialization stage that is also of $\mathcal{O}(N)$. This method makes use of the analytical and low-rank property of $\boldsymbol{A}$ as well as hierarchical matrices, and it is the method that we will attempt to improve in Sec.~\ref{sec:faster}. Keiner \cite{keiner09} extends this method to transforms between families of Gegenbauer polynomials and improves the error estimates of the approximations. Shen et al. \cite{shen19} describe a similar method for Jacobi polynomials, also based on the low-rank property of the connection matrices and using hierarchical semiseparable matrices. 
The algorithm of Shen et al. also requires $\mathcal{O}(N)$ operations for evaluation after a longer $\mathcal{O}(N^2)$ initialization stage. Keiner \cite{keiner11} also makes use of semiseparable matrices and describes an $\mathcal{O}(N \log N)$ divide and conquer approach. Hale and Townsend \cite{Hale2014} describe a $\mathcal{O}(N(\log N)^2 / \log \log N)$ method based on Stieltjes' asymptotic formula for Legendre polynomials \cite{stieltjes}.  This method is sold on the merit of requiring no initialization (except from the FFTs involved) and being simple enough to be implemented in 100 lines of MATLAB code. An even simpler $\mathcal{O}(N \log^2 N)$ approach is described by Townsend et al.~\cite{Townsend18}, decomposing the connection matrix into diagonally scaled Hadamard products involving Toeplitz and Hankel matrices. This method is implemented by the MATLAB software system Chebfun \cite{chebfun} and the Julia package called FastTransforms.jl \cite{fasttransforms}. The method to be described in Sec.~\ref{sec:faster} will be shown to be significantly faster and can also be implemented in approximately 100 lines of compact Python code. 

\section{A direct method}
\label{sec:direct}

The structure of $\boldsymbol{A}$ (see Eqs. \eqref{eq:hik} and \eqref{eq:matH}) allows for a low memory and quick direct approach for computing $\boldsymbol{Af}^{leg}$ and thus $\boldsymbol{f}^{cheb}$ from $\boldsymbol{f}^{leg}$. The method applies diagonalwise as shown in Alg. \ref{alg:l2cdirect}. The method is very efficient for small $N$, but naturally it will become expensive for large $N$ since the number of floating point operations (flops) of a direct matrix vector product scales like $\mathcal{O}(N^2)$. 
We note that a direct method that precomputes all $\lambda_{n/2}\lambda_{n/2+i}$ (see Alg. \ref{alg:l2cdirect}) will lead to fewer floating point operations for the execution, but it is not necessarily faster, because the alternative leads to fewer memory access operations. To see this consider what the matrix $\boldsymbol{A}$ looks like in Eq. \eqref{eq:matH}:
\begin{equation}
\boldsymbol{A} =
\begin{bmatrix}
\lambda_0\lambda_0 & 0      & \lambda_1 \lambda_1 & 0      & \lambda_2\lambda_2 & 0      &  \cdots \\
0      & \lambda_0\lambda_1 & 0       & \lambda_1\lambda_2 & 0      & \lambda_2\lambda_3  &\cdots\\
0      & 0      & \lambda_0 \lambda_2 & 0      & \lambda_1\lambda_3 & 0      & \cdots \\
0      & 0      & 0       & \lambda_0\lambda_3 & 0      & \lambda_1\lambda_4 & \cdots\\
0      & 0      & 0       & 0      & \lambda_0\lambda_4 & 0      & \cdots \\
0      & 0      & 0       & 0      & 0      & \lambda_0 \lambda_5 & \cdots\\
\vdots & \vdots & \vdots  & \vdots & \vdots & \vdots & \ddots
\end{bmatrix}. \label{eq:matH}
\end{equation}
In order to compute with the 12 nonzero items of $\boldsymbol{A}$ shown, it is sufficient for the compiler to retrieve only 6 numbers from memory $(\lambda_0, \lambda_1, \ldots, \lambda_5)$. So, compared to precomputing all matrix items, memory access for the matrix is halved (for these items), whereas the flop count is increased from two (one multiplication and one addition) to three (two multiplications and one addition) per non-zero item. Whichever is faster needs to be tested on a given computing system.

\begin{algorithm}[h!]
\caption{A direct Legendre-to-Chebyshev transform 
$\boldsymbol{f}^{cheb} = \boldsymbol{C}^{-1} \boldsymbol{A} \boldsymbol{f}^{leg}$ }
\label{alg:l2cdirect}
\begin{algorithmic}[1]
\Function{DIRECTL2C}{$\boldsymbol{f}^{leg},\boldsymbol{\lambda},\boldsymbol{C}$}
\Input{$\boldsymbol{f}^{leg}$}{array $\in \mathbb{R}^{N}$}{Legendre coefficients}
\Input{$\boldsymbol{\lambda}$}{array $\in \mathbb{R}^{N}$}{$\lambda_k = \Lambda(k)$, see Eq. \eqref{eq:Lambda}}
\Input{$\boldsymbol{C}$}{matrix $\in \mathbb{R}^{N \times N}$}{diagonal matrix, see Eq. \eqref{eq:AfCf}}
\Output{$\boldsymbol{f}^{cheb}$}{array $\in \mathbb{R}^{N}$}{Chebyshev coefficients}
\State $N \gets \mathbf{{len}}~\boldsymbol{f}^{leg}$
\State ${f}^{cheb} \gets 0$
\State $n \gets 0$
\While {$n < N$}
    \For{$i \gets 0, N-n-1$}
        \State ${f}^{cheb}_i \gets {f}_i^{cheb} + \lambda_{n/2}\lambda_{n/2+i} {f}^{leg}_{i+n}$
    \EndFor
    \State $n \gets n+2$
\EndWhile
\State $ \boldsymbol{f}^{cheb} \gets \boldsymbol{C}^{-1} \boldsymbol{f}^{cheb}$
\State \Return $\boldsymbol{f}^{cheb}$
\EndFunction
\end{algorithmic}
\end{algorithm}


\section{A Fast Multipole Method}
\label{sec:FMM}
We will first describe a fast $\mathcal{O}(N \log_2 N)$ multipole method for the matrix vector product $\boldsymbol{A}\boldsymbol{f}^{leg}$ that is very easy to implement.\footnote{The method for the reverse direction is basically identical and not shown in detail.} The method uses the low-rank properties of $\boldsymbol{A}$ and the general idea from \cite{alpert91}. However, where \cite{alpert91} approximates $\boldsymbol{A}$ using nodal Lagrange polynomials in real space, we will take a modal approach and make use of Chebyshev polynomials in spectral space. In an effort to simplify notation, we will in this section write the matrix vector product as
\begin{equation}
    \boldsymbol{z} = \boldsymbol{A} \boldsymbol{f}. \label{eq:matvecsimpl}
\end{equation}

\begin{figure}
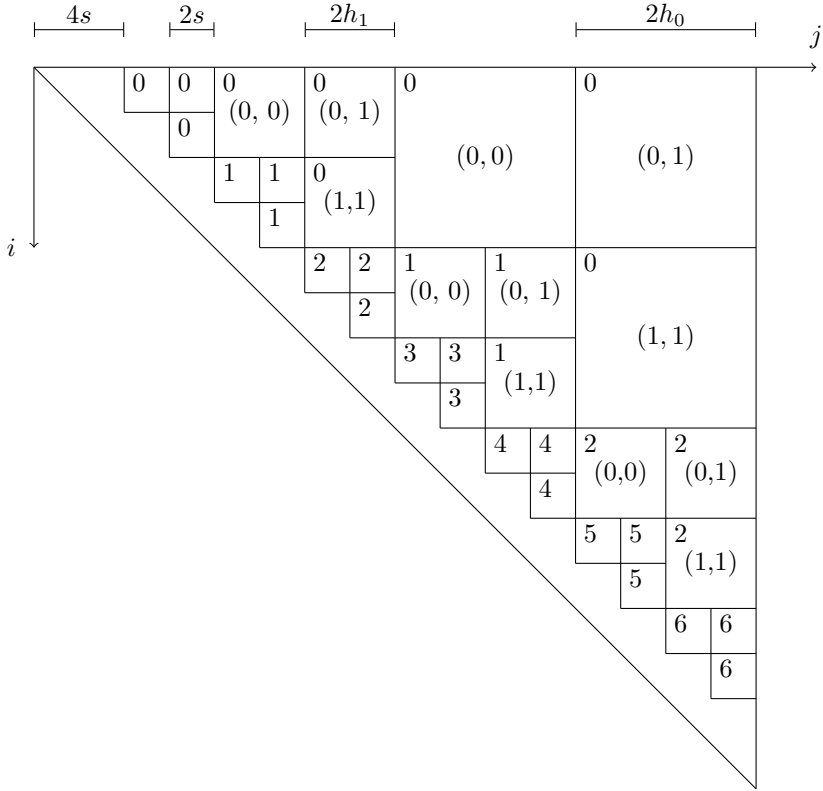

    \centering
    \ctikzfig{levels2}
    \caption{Example of a decomposition of the matrix $\boldsymbol{A}$ using three levels. Level 0 is the one with the fewest blocks, and then levels 1 and 2 have gradually more. For levels 0 and 1 the local block numbers are given in the upper left corner of each square, whereas the local indices for a submatrix on a block are shown in parenthesis. For level 2 only the local block number is printed.}
    \label{fig:3levels}
\end{figure}

\subsection{Decomposition into levels, blocks and square submatrices}
\label{sec:decomposition}
The main idea of the fast multipole method is to decompose the matrix $\boldsymbol{A}$ into smaller submatrices, and then use low-rank approximations to efficiently compute the matrix-vector product of each submatrix.
For decomposition of $\boldsymbol{A}$ we will use a system based on \emph{levels}, \emph{blocks} and square \emph{submatrices}. Figure \ref{fig:3levels} shows an example where the upper triangular part of $\boldsymbol{A}$ has been decomposed using three levels, i.e., there are square submatrices of exactly three different sizes. The largest squares represent level 0, and the two gradually smaller are thus levels 1 and 2. We will use $\gamma \in \mathbb{N}$ to refer to a specific level and the total number of levels is $L \in \mathbb{N}$. The structured decomposition of the matrix into levels, blocks and submatrices will allow us to easily create iterators that runs over all blocks, submatrices and vectors on any or all levels. Note that the unmarked part of the matrix closest to the main diagonal in Fig. \ref{fig:3levels} will be treated with a direct approach, as described in Sec~\ref{sec:direct}.

On any given level we define a \emph{block} as the three neighbouring \emph{submatrices} of the same shape that share at least one common edge with at least one other submatrix in the block. The local (to each level) block numbers are given in Fig. \ref{fig:3levels} in the upper left corner of each square, and for level 2 the local block numbers are the only numbers printed.
There are 11 blocks altogether in Fig. \ref{fig:3levels}, as we are counting 1, 3 and 7 blocks on levels 0, 1 and 2, respectively. 

All square submatrices on level ${\gamma}$ have shape $2h_{\gamma} \times 2h_{\gamma}$, as shown in Fig. \ref{fig:3levels}, where the even shape $2h_{\gamma}$ has been chosen since the matrix will be further decomposed into equally shaped odd and even parts. The size of the squares on each level, $h_{\gamma}$, can be computed as
\begin{equation}
    h_{\gamma} = s \, 2^{L-\gamma-1},
\end{equation}
where $s=h_{L-1}$ is half the size of the smallest matrices, which needs to be specified.

The numbers in parenthesis in Fig. \ref{fig:3levels}, in the center of the squares on levels 0 and 1, show the local indices $(p,q) \in \mathbb{N}^2$ used for a square \emph{submatrix} on a block. For efficient storage and lookup these two indices will be mapped into a single index $r$ using a column-major numbering scheme
\begin{equation}
r(p, q) = p + q(q+1)/2, \quad  r \in \{0, 1, 2\} \text{ and } q \ge p,     \label{eq:cfrompq}
\end{equation}
such that $r=0, 1 \text{ and } 2 \implies (p, q)=(0, 0), (0, 1) \text{ and } (1, 1)$, respectively.

We access a submatrix using the level number ${\gamma}$, local block number $b$ and the local submatrix index $r(p,q)$
\begin{equation}
    \boldsymbol{{A}}({\gamma},b,r) = (a_{i+m, j+n})_{m,n=0}^{2h_{\gamma}-1} \in \mathbb{R}^{2h_{\gamma} \times 2h_{\gamma}}, \label{eq:Aunderline}
\end{equation}
where the integers $i$ and $j$ represent the global indices of the upper left corner of the submatrix $\boldsymbol{{A}}(\gamma, b, r)$.
It is easily shown that these global indices can be computed as
\begin{equation}
        i(\gamma, b, p) = 2h_{\gamma}(2b + p) \in \mathbb{N} \, \text{ and } \,
        j(\gamma, b, q) = 2h_{\gamma}(2b+q+2) \in \mathbb{N}. 
    \label{eq:globalij}
\end{equation}
Similarly, the odd ($\sigma=1$) and even ($\sigma=0$) parities of the same submatrix will be accessed as
\begin{equation}
    \boldsymbol{{A}}^{\sigma}({\gamma},b,r) = ({a}_{i+2m+\sigma, j+2n+\sigma})_{m,n=0}^{h_{\gamma}-1} \in \mathbb{R}^{h_{\gamma} \times h_{\gamma}}.\label{eq:Aoverline}
\end{equation}

The implementation of the FMM is based on an iterator that runs over all levels, blocks and submatrices in an ordered fashion. In order to create such an iterator we simply need to know the number of levels $L$ and the number of blocks on each level $b_{\gamma}$. 
The number of blocks on a given level $\gamma$ can be computed as twice the number of blocks on level $\gamma-1$ plus one. Starting from $b_0=1$, this becomes through simple recursion
\begin{equation}
    b_{\gamma} = 2^{\gamma+1}-1. \label{eq:Nb}
\end{equation}
Note that the chosen decomposition restricts the possible shapes $N$ of the input vector $\boldsymbol{f}\in \mathbb{R}^N$, because the algorithm requires that
\begin{equation}
    N = 2s(1+b_{L}) = s 2^{L+2}. \label{eq:Ntot}
\end{equation}
We may use Eq. \eqref{eq:Ntot} to compute the number of levels required by any input vector $\boldsymbol{f}$ of length $\lvert \boldsymbol{f}\rvert$ not necessarily satisfying \eqref{eq:Ntot}
\begin{equation}
    L = \lceil \log_2 \frac{\lvert \boldsymbol{f} \rvert}{\hat{s}} \rceil -2. \label{eq:numlevels}
\end{equation}
Here $\hat{s}$ is a chosen (optimal) size of the smallest submatrices and $\lceil x \rceil$ represents the ceiling function mapping $x$ to the nearest integer greater than or equal to $x$. We then recompute $s = \lceil \lvert \boldsymbol{f} \rvert /2^{L+2} \rceil  \in [\sfrac{\hat{s}}{2}, \hat{s}]$ and finally set $N = s 2 ^{L+2}$, where $N \ge \lvert \boldsymbol{f} \rvert$.
The bottom line is that any input array $\boldsymbol{f}$ needs to be padded with $s2^{L+2}-\lvert\boldsymbol{f}\rvert$ zeros for the algorithm to work. In what follows we will simply assume that $N$ satisfies Eq.~\eqref{eq:Ntot} and thus $L=\log_2\frac{N}{s}-2$.

The total number of blocks $N_b$ and submatrices $N_c$ are obtained by summing $b_{\gamma}$ over all levels, and we obtain
\begin{equation}
    N_b = \frac{N}{2s}-L-2 \quad \text{and} \quad N_c = 3N_b, \label{eq:total_blocks_and_mats}
\end{equation}
respectively. Since the number of levels $L$ grows as $\mathcal{O}(\log_2 N)$, we see that both the number of blocks and submatrices in Eq. \eqref{eq:total_blocks_and_mats} grow asymptotically as $\mathcal{O}(N)$.


The global vectors $\boldsymbol{z} \in \mathbb{R}^N$ and $\boldsymbol{f} \in \mathbb{R}^N$, that take part in the matrix vector product \eqref{eq:matvecsimpl}, are also looked up locally using the global indices $i$ and $j$ from Eq. \eqref{eq:globalij}. As such, we can easily find the parts of the global vectors that belong to any given level $\gamma$, block $b$ and submatrix $r(p,q)$
\begin{equation}
\boldsymbol{f}(\gamma, b, q) = \{f_{j+m}\}_{m=0}^{2h_{\gamma}-1} \in \mathbb{R}^{2h_{\gamma}}, \,
\boldsymbol{z}(\gamma, b, p) = \{z_{i+m}\}_{m=0}^{2h_{\gamma}-1} \in \mathbb{R}^{2h_{\gamma}}.
\end{equation}
Likewise, we get the larger vector blocks associated with entire blocks or levels as 
\begin{align}
\boldsymbol{f}(\gamma, b) &= \{f_{j_b+m}\}_{m=0}^{4h_{\gamma}-1} \in \mathbb{R}^{2 \times 2h_{\gamma}}, \, \boldsymbol{f}(\gamma) = \{f_{4h_{\gamma}+m}\}_{m=0}^{N-4h_{\gamma}-1} \in \mathbb{R}^{b_{\gamma} \times 2 \times 2h_{\gamma}}, \label{eq:fgamma} \\
\boldsymbol{z}(\gamma, b) &= \{z_{i_b+m}\}_{m=0}^{4h_{\gamma}-1} \in \mathbb{R}^{2 \times 2h_{\gamma}}, \, \boldsymbol{z}(\gamma) = \{z_{m}\}_{m=0}^{N-4h_{\gamma}-1} \in \mathbb{R}^{b_{\gamma} \times 2 \times 2h_{\gamma}}, \label{eq:zgamma}
\end{align}
where $i_b(b, \gamma)=4bh_{\gamma}$ and $j_b(b, \gamma)=4(1+b)h_{\gamma}$ represent the indices of the upper left corner of block $b$ on level $\gamma$.
We treat the arrays in \eqref{eq:fgamma} and \eqref{eq:zgamma} as row-major, multidimensional arrays simply for convenience in future algorithms. However, the underlying data are just a single, contiguous block of numbers, and we might as well have been using, e.g., $ \boldsymbol{f}(\gamma) \in \mathbb{R}^{4b_{\gamma}h_{\gamma}}$.
For all vector blocks, if the vector is also endowed with a $\sigma$ superscript, then a subdivision is further made into odd and even parts, like $\boldsymbol{f}^{\sigma}( \gamma, b, q) = \{f_{j+2m+\sigma}\}_{m=0}^{h_{\gamma}-1} \in \mathbb{R}^{h_{\gamma}}$.

\begin{figure}
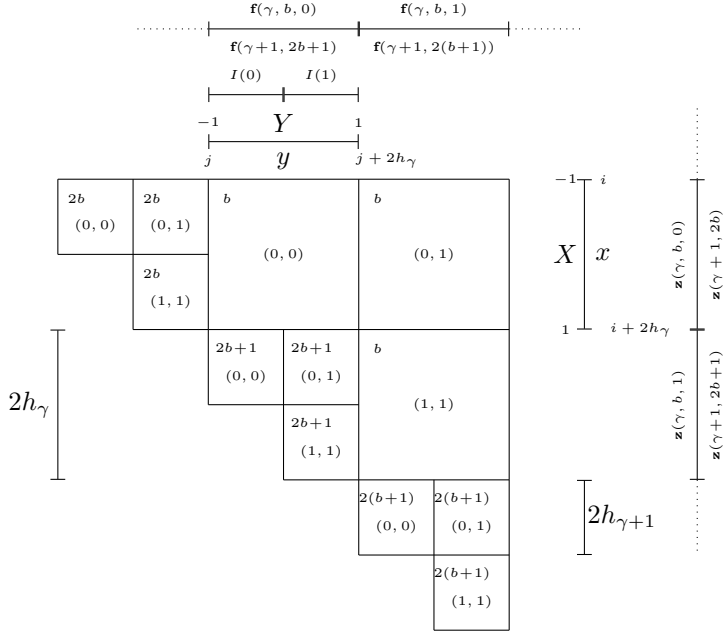

    \centering
    \ctikzfig{2levelsFMM}
    \caption{Example of a local part of the matrix $\boldsymbol{A}$, showing two neighbouring levels. There is one block on level $\gamma$, with block number $b$. The block numbers on level $\gamma+1$ are shown in upper left corners as computed from $b$. Also shown are the reference domains $X$ and $Y$ on level $\gamma$. For submatrix $(0, 0)$ on level $\gamma$, where the upper left corner is at global indices $i, j$, we also show the subintervals $I(0)$ and $I(1)$ used in Eq. \eqref{eq:wtk}.} 
    \label{fig:2levelsFMM}
\end{figure}

Figure \ref{fig:2levelsFMM} gives a more detailed and generic description of a zoomed-in part of a larger global matrix $\boldsymbol{A}$, illuminating how the global vectors on different levels are related. A correlation between vector blocks on neighbouring levels is obviously
\begin{align}
    \boldsymbol{f}(\gamma, b, q) &= \boldsymbol{f}(\gamma+1, 2b+q+1), \label{eq:flevels} \\ 
    \boldsymbol{z}(\gamma, b, p) &= \boldsymbol{z}(\gamma+1, 2b+p), \label{eq:zlevels}
\end{align}
where the equalities should be understood in terms of the underlying data.

\subsection{The 2D forward Chebyshev transform}
The most important part of the FMM is to approximate a submatrix $\boldsymbol{A}(\gamma, b, r)$ with a Chebyshev series. As described by Alpert and Rokhlin \cite{alpert91}, we may view the global matrix $\boldsymbol{A}$ as a continuous function $\mathcal{A}(x, y): \mathbb{R} \times \mathbb{R} \rightarrow \mathbb{R}$, such that $\mathcal{A}(i, j) = a_{ij}$ for all nonzero items of $\boldsymbol{A}$. The continuous function 
\begin{equation}
    \mathcal{A}(x, y) = \Lambda\left(\frac{y-x}{2}\right) \Lambda\left(\frac{y+x}{2}\right),
\end{equation}
is smooth for $y>x$. 

For all square submatrices $\boldsymbol{{A}}(\gamma, b, r)$ we will approximate  $\mathcal{A}(x, y)$ in the corresponding domain $\Omega(\gamma,b,r) = [i, i+2h_{\gamma}] \times [j, j+2h_{\gamma}]$, by projection to a two-dimensional reference tensor product space $W = V_M(I) \otimes V_M(I)$, where $I = [-1, 1]$, $V_M(I) = \text{span}\{T_k(\MakeUppercase{x})\}_{k=0}^{M-1}$ and $\MakeUppercase{x} \in I$ is a reference coordinate. For any square submatrix with upper left corner in $(i, j)$, there is an affine mapping from the physical coordinates $(x, y) \in \Omega$ to the reference coordinates 
\begin{equation}
 \MakeUppercase{x}=-1+(x-i)/h_{\gamma} \, \text{ and } \,
 \MakeUppercase{y}=-1+(y-j)/h_{\gamma}, \label{eq:Xfromx}
\end{equation} 
where $(i,j)$ as before are given by Eq. \eqref{eq:globalij} and $Y$ is used to identify a reference coordinate for the second axis.

An expansion for $\mathcal{A}(x,y)$ in the two-dimensional Chebyshev space $W$ is 
\begin{equation}
    \mathcal{A}(x, y) \approx \sum_{k=0}^{M-1} \sum_{l=0}^{M-1} \hat{a}_{kl} T_{k}(\MakeUppercase{x})T_{l}(\MakeUppercase{y}), \quad (\MakeUppercase{x}, \MakeUppercase{y}) \in I^2. \label{eq:Aij}
\end{equation}
The coefficients $\boldsymbol{\hat{A}} = (\hat{a}_{kl})_{k,l=0}^{M-1} \in \mathbb{R}^{M \times M}$ are found by $L^2_{\omega}(I^2)$ orthogonal projection, or, more simply put, from a 2D forward Chebyshev transform. Multiplying Eq. \eqref{eq:Aij} through by a test function $T_m(X)T_n(Y)$ and weights $\omega_X = (1-X^2)^{-\sfrac{1}{2}}$ and $ \omega_Y = (1-Y^2)^{-\sfrac{1}{2}}$, integrating over the domain and using the orthogonality of Chebyshev polynomials leads to
\begin{equation}
    \hat{a}_{kl} = \frac{4}{\pi^2 c_k c_l} \left( \left({\mathcal{A}}(x, y), T_k(\MakeUppercase{x})\right)_{\omega_X}, T_l(\MakeUppercase{y}) \right)_{\omega_Y}, \, \forall \, (k, l) \in \mathcal{I}^2_M, \label{eq:akl2}
\end{equation}
where $\mathcal{I}_M=\{0, 1, \ldots, M-1\}$. 

Keiner \cite{keiner09} shows that it is sufficient to use $M=18$ for full double precision for the Chebyshev approximation of $\mathcal{A}(x, y)$ on any of the problems square submatrices. Hence the scalar products can also be computed without any loss of accuracy using numerical quadrature (instead of exact integration) on Chebyshev-Gauss points $ \boldsymbol{X}^G = \boldsymbol{Y}^G = \{\cos \left((m+\sfrac{1}{2})\pi/M\right)\}_{m=0}^{M-1}$, with constant quadrature weights $\{\pi/M\}_{m=0}^{M-1}$, and the definition $T_k(x) = \cos(k \cos^{-1}x)$, such that Eq. \eqref{eq:akl2} becomes
\begin{equation}
    \hat{a}_{kl} = \frac{4 }{c_kc_l M^2} \sum_{m=0}^{M-1}\sum_{n=0}^{M-1}\mathcal{A}(x_m, y_n) \cos(k(m+\sfrac{1}{2})\pi/M) \cos(l(n+\sfrac{1}{2})\pi/M), \label{eq:aklmore}
\end{equation}
which is easily computed with a discrete cosine transform of type 2. Note that $x_m = i + h_{\gamma} (X^G_m + 1)$ and $y_n = j + h_{\gamma}(Y^G_n+1)$ from Eq. \eqref{eq:Xfromx}.

We find the coefficients for all submatrices in $\boldsymbol{A}$ by running over all levels, blocks and submatrices, evaluating the matrix $\left(\mathcal{A}(x_m, y_n)\right)\in \mathbb{R}^{M\times M}$ and performing the discrete cosine transform for each square. This is the major cost of the initialization stage of the FMM. Since the initialization cost for each submatrix is the same (without using any optimizations), the total initialization cost will scale as $\mathcal{O}(N)$ since the number of submatrices grows asymptotically as $\mathcal{O}(N)$, see Eq. \eqref{eq:total_blocks_and_mats}. In Sec. \ref{sec:fastinit} we will discuss some significant optimizations for this initialization.

Note that the original multipole method of Alpert et al. \cite{alpert91} used nodal Lagrange polynomials $\{\ell_k(\MakeUppercase{x})\}_{k=0}^{M-1}$ instead of Chebyshev,\footnote{To be perfectly consistent they also used the reference interval $[0, 1]$ and not $[-1, 1]$.} and the approximation
\begin{equation}
    \mathcal{A}(x, y) \approx \sum_{k=0}^{M-1} \sum_{l=0}^{M-1} \mathcal{A}(x_k, y_l) \ell_{k}(\MakeUppercase{x}) \ell_{l}(\MakeUppercase{y}), \label{eq:AijLag}
\end{equation}
instead of Eq. \eqref{eq:Aij}. Hence, for the original method the initializtion cost is smaller, since it does not require the Chebyshev transform of the otherwise same matrix $(\mathcal{A}(x_m, y_n))_{m,n=0}^{M-1}$.


\subsection{Evaluating the matrix vector product on one submatrix}
\label{sec:submatrix1}

For any given submatrix $\boldsymbol{\underline{A}}=\boldsymbol{{A}}(\gamma,b,r)\in \mathbb{R}^{2h_{\gamma} \times 2h_{\gamma}}$ the matrix vector product that we are interested in can be computed as 
\begin{equation}
    \underline{z}_m = \sum_{n=0}^{2h_{\gamma}-1} \underline{a}_{mn} \underline{f}_n, \quad m=0,1, \ldots, 2h_{\gamma}-1, \label{eq:matvec0}
\end{equation}
where $\boldsymbol{\underline{z}} = \boldsymbol{z}(\gamma, b, p) $, $\boldsymbol{\underline{f}}=\boldsymbol{{f}}(\gamma, b, q)$ and an underline notation is used in general to indicate a sub-vector/matrix.

If we do not take the sparsity of $\boldsymbol{\underline{A}}$ into consideration the flop count of \eqref{eq:matvec0} will be $8h_{\gamma}^2$. Note that we do not attempt to find an exact flop count and will only count the flops of the most deeply nested statements, using Table 1.1.2 of \cite{GoVa13} for reference.  
For efficiency we now also split into odd and even parities:
\begin{equation}
    \underline{{z}}^{\sigma}_m = \sum_{n=0}^{h_{\gamma}-1} \underline{{a}}^{\sigma}_{mn} \underline{{f}}^{\sigma}_n, \quad \forall \, m \in \mathcal{I}_{h_{\gamma}} \text{ and } \sigma \in \{0, 1\},
\end{equation}
where $\mathcal{I}_{h_{\gamma}} = \{0, 1, \ldots, h_{\gamma}-1\}$, $\underline{z}_{m}^{\sigma} = \underline{z}_{2m+\sigma},  \underline{f}_{n}^{\sigma} = \underline{f}_{2n+\sigma}$ and $\underline{a}_{mn}^{\sigma} = \underline{a}_{2m+\sigma, 2n+\sigma}$. Since all zero items of $\boldsymbol{A}$ now have been eliminated, the cost has been reduced to $4h_{\gamma}^2$. 

We now come to the actual multipole method and replace $\underline{a}^{\sigma}_{mn}$ by its Chebyshev approximation \eqref{eq:Aij}, such that the matrix vector product becomes
\begin{equation}
    \underline{z}^{\sigma}_m = \sum_{n=0}^{h_{\gamma}-1} \sum_{k=0}^{M-1} \sum_{l=0}^{M-1} \hat{a}_{kl} T_{k}(X^{\sigma, \gamma}_m)T_{l}(Y^{\sigma, \gamma}_n)  \underline{f}^{\sigma}_n, \quad m \in \mathcal{I}_{h_{\gamma}}, \sigma \in \{0, 1\}, \label{eq:matveczs}
\end{equation}
where the Vandermonde matrix $\boldsymbol{T}(\sigma, \gamma) = \left(T_k(X^{\sigma,\gamma}_m)\right)_{m, k = 0}^{h_{\gamma}-1, M-1} \in \mathbb{R}^{h_{\gamma} \times M}$, where 
\begin{equation}
    X^{\sigma,\gamma}_m=Y^{\sigma,\gamma}_m=-1+(2m+\sigma)/h_{\gamma},\quad m \in \mathcal{I}_{h_{\gamma}},
    \label{eq:Xs}
\end{equation} 
is the (uniform) discretization of the local reference coordinate. Note that the Chebyshev coefficients $\hat{a}_{kl}$ are the same for the odd or even computations, and the parity only affects the Vandermonde matrix. Also note that there is only one Vandermonde matrix for a given parity since the matrices $(T_{k}(X^{\sigma,\gamma}_m))$ and $(T_{l}(Y^{\sigma,\gamma}_n))$ in \eqref{eq:matveczs} are identical. 

If we now rewrite \eqref{eq:matveczs} with matrix notation, $\boldsymbol{\hat{A}}=\boldsymbol{\hat{A}}(\gamma, b, r)$ and $\boldsymbol{T} = \boldsymbol{T}(\sigma, \gamma)$: 
\begin{equation}
    \boldsymbol{\underline{z}}^{\sigma} = \boldsymbol{{T}}\boldsymbol{\hat{A}} \boldsymbol{{T}}^T \boldsymbol{\underline{f}}^{\sigma}, \quad \forall \, {\sigma}\in \{0, 1\},
\end{equation}
we see that for given ${\sigma}$ the cost is $2Mh_{\gamma}$ flops to do $ \boldsymbol{w} = \boldsymbol{{T}}^T \boldsymbol{\underline{f}}^{\sigma}$, where $\boldsymbol{w}  \in \mathbb{R}^{M}$ is a work array, $2M^2$ to do $ \boldsymbol{c} = \boldsymbol{\hat{A}} \boldsymbol{w}$ for another work array $\boldsymbol{c} \in \mathbb{R}^M$, and $2h_{\gamma}M$ to do $\boldsymbol{\underline{z}}^{\sigma} = \boldsymbol{{T}} \boldsymbol{c}$, for a total cost of $4Mh_{\gamma}+2M^2$. Adding the two odd and even parities together, the cost becomes $8Mh_{\gamma}+4M^2$, so if $M << h_{\gamma}$ there can be a significant speedup over the $4h_{\gamma}^2$ cost of the direct computation. 

\subsection{Evaluating the matrix vector product for a block of submatrices}
\label{sec:submatrix2}

For any given \emph{block} $b$ of submatrices the matrix vector product that we are interested in can be computed by running over all 3 submatrices in the block and adding up the contributions of each, as shown in Alg. \ref{alg:matvecblock}. The cost of this operation is less than 3 times $8Mh_{\gamma}+4M^2$ (which is the cost for one submatrix isolated), because the vector block $\boldsymbol{\underline{f}}(q)$ is independent of $p$ and $\boldsymbol{\underline{z}}(p)$ is independent of $q$. We get the flop count as detailed in Lemma \ref{lem:blockflops}.
\begin{algorithm}[h!]
\caption{Matrix vector product for a single block of submatrices for a given parity.}
\label{alg:matvecblock}
\begin{algorithmic}
\Function{MATVECBLOCK}{$\boldsymbol{f}, \boldsymbol{T}, \boldsymbol{\hat{A}}, \boldsymbol{z}$}
\Input{$\boldsymbol{{f}}$}{$\mathbb{R}^{2 \times h}$}{Legendre coefficients for given block and parity}
\Input{$\boldsymbol{{T}}$}{$\mathbb{R}^{h \times M}$}{Vandermonde matrix for given level and parity}
\Input{$\boldsymbol{\hat{A}}$}{$3$ sequence of $\mathbb{R}^{M \times M}$}{Chebyshev matrices of coefficients}
\Output{$\boldsymbol{{z}}$}{$\mathbb{R}^{2 \times h}$}{Chebyshev coefficients for given block and parity}
\State $\boldsymbol{{c}} \gets 0$ \Comment{Initialize work array $\boldsymbol{c}\in \mathbb{R}^{2 \times M}$}
\For{$q \gets 0, 1$}
    \State $\boldsymbol{w} \gets \boldsymbol{{T}}^T \boldsymbol{{f}}(q)$ \Comment{Work array $\boldsymbol{w} \in \mathbb{R}^{M}$}
    \For{$p \gets 0, q$}
        \State $r \gets p+q(q+1)/2$ \Comment{Eq. \eqref{eq:cfrompq}}
        \State $\boldsymbol{{c}}(p) \gets  \boldsymbol{{c}}(p) + \boldsymbol{\hat{A}}(r) \boldsymbol{w}$
    \EndFor
\EndFor
\For{$ p \gets 0, 1$}
    \State $\boldsymbol{{z}}(p) \gets \boldsymbol{{z}}(p) + \boldsymbol{{T}} \boldsymbol{{c}}(p)$
\EndFor
\State \Return $\boldsymbol{{z}}$
\EndFunction
\end{algorithmic}
\end{algorithm}

\begin{lemma}
\label{lem:blockflops}
The significant flop count for a matrix vector product over a block is
\begin{equation}
    16 h_{\gamma}M + 12M^2. \label{eq:numflops}
\end{equation}

\begin{proof}
    We can count the flops line by line in Alg. \ref{alg:matvecblock}, taking only the most deeply nested computations into account. The first operation ${\boldsymbol{{T}}}^T \boldsymbol{{f}}(q)$ requires $8 h_{\gamma}M$ flops since the matrix vector product costs $2h_{\gamma}M$ and it is computed $2$ times for each parity. The second step $\boldsymbol{{c}}(p) \gets \boldsymbol{{c}}(p) + \boldsymbol{\hat{A}}(r(p, q)) \boldsymbol{{w}}$ is a matrix vector product of cost $2M^2$, computed $3$ times per parity. The third and final step $ {\boldsymbol{{T}}}\boldsymbol{{c}}(p)$  costs $8h_{\gamma}M$ flops for the same reasons as the first step. Counting all three steps yields the number in Eq. \eqref{eq:numflops}.
\end{proof}
\end{lemma}

\begin{algorithm}
\caption{Procedure used to compute matrix vector product $\boldsymbol{z} = \boldsymbol{A} \boldsymbol{f}$.}    \label{alg:complete}
\begin{algorithmic}
    \For{$\sigma \gets 0, 1$}
        \For{$\gamma \gets 0, L-1$}
            \For{$b \gets 0, b_{\gamma}-1$}
                \State $\boldsymbol{z}^{\sigma}(\gamma, b) \gets $ MATVECBLOC($\boldsymbol{f}^{\sigma}(\gamma, b)$, 
                $\boldsymbol{T}(\sigma, \gamma)$, 
                ${\boldsymbol{\hat{A}}}(\gamma, b)$,
                $\boldsymbol{z}^{\sigma}(\gamma, b)$)
            \EndFor
        \EndFor
    \EndFor
\end{algorithmic}
\end{algorithm}

The complete matrix vector product can now be computed simply by running over all levels and blocks for odd and even parities of the matrix, as shown in the simplified Alg. \ref{alg:complete}.
The total cost of the method can then be estimated as
\begin{equation}
    \mu s N + \sum_{\gamma=0}^{L-1} b_{\gamma} \left(16 h_{\gamma}M + 12M^2 \right), \label{eq:totalnumflops}
\end{equation}
where $\mu s N$ represents the cost of using the direct method for the part of the matrix that is close to the diagonal and separated from the submatrices. The scaling factor $\mu=3$ if the elements of the matrix $\boldsymbol{A}$ are precomputed. However, $\mu=4.5$ if we simply use the vector $\boldsymbol{\lambda}$ as in Alg. \ref{alg:l2cdirect}. By inserting for $b_{\gamma}$ and $h_{\gamma}$ we can simplify Eq. \eqref{eq:totalnumflops} further and obtain a total cost lower than
\begin{equation}
    (\mu s+4LM+\frac{6M^2}{s})N, \label{eq:totalnumflops2}
\end{equation}
where the last number represents $4N_c M^2 = 12(\frac{N}{2s}-L-2)M^2 < 6NM^2/s$, or the cost of applying all $N_c$ submatrices $\boldsymbol{\hat{A}}(\gamma, b, r)$, each with a cost of $2M^2$ per parity. Since the number of levels $L \propto \log_2 N/s$, we get that the total cost of this very simple method scales as $\mathcal{O}(N \log_2 N)$.


\section{A Faster Multipole Method}
\label{sec:faster}
The Fast Multipole Method described in Sec. \ref{sec:FMM} is not optimal since the relatively expensive matrix vector products $ {\boldsymbol{{T}}}^T \boldsymbol{\underline{f}}^{\sigma}(q)$ and ${\boldsymbol{{T}}} \boldsymbol{{c}}(p)$ in Alg. \ref{alg:matvecblock} are computed at each level. In the method described by Alpert and Rokhlin \cite{alpert91} this obstacle is overcome by their Eq.~(4.8), which maps the Lagrange interpolating polynomials on one level to those on the next finer level. Below we will adapt a similar approach for the Chebyshev polynomials and this modification will make the total method scale as $\mathcal{O}(N)$. The complete faster algorithm is given in Alg. \ref{alg:fastercomplete} and it is described in some more detail below.

\subsection{The algorithm}
The first objective of the faster method is to speed up the matrix vector product $ \boldsymbol{w} = {\boldsymbol{{T}}}^T \boldsymbol{\underline{f}}^{\sigma}(q)$, which is computed $2$ times per block on all levels. We thus start by predefining a multidimensional work array to hold this result for all levels: $\boldsymbol{w}(\gamma) \in \mathbb{R}^{b_{\gamma} \times 2 \times M} \, \forall \, \gamma \in \mathcal{I}_L$. 
For any level, block and parity the work vector $\boldsymbol{w}(\gamma, b, q) \in \mathbb{R}^M$ can be computed directly as the matrix-vector product
\begin{equation}
    \boldsymbol{w}(\gamma, b, q) = {\boldsymbol{T}}^T(\sigma, \gamma) \boldsymbol{f}^{\sigma}(\gamma, b, q) \quad \in \mathbb{R}^M, \label{eq:wTf}
\end{equation}
costing $2h_{\gamma}M$ operations per parity.
As for the global vectors $\boldsymbol{f}$ and $\boldsymbol{z}$,  we here look up the local part of the work vector using level $\gamma$, block $b$ and local column index $q$. Since it is merely a work vector, we drop the parity superscript on $\boldsymbol{w}$.

The direct approach for computing $\boldsymbol{w}(\gamma)$ (step (1) in Alg. \ref{alg:fastercomplete}) will only be used for the highest level $\gamma=L-1$, where it can be set up efficiently as a tensor contraction
\begin{equation}
    \boldsymbol{w}(L-1) = \boldsymbol{f}^{\sigma}(L-1) \boldsymbol{T}(\sigma, L-1)  \quad \in \mathbb{R}^{b_{L-1} \times 2 \times M}, \label{eq:wMM}
\end{equation}
between the two tensors $\boldsymbol{f}^{\sigma}(L-1) \in \mathbb{R}^{b_{L-1} \times 2 \times s}$ and $\boldsymbol{T}(\sigma, L-1)\in \mathbb{R}^{s \times M}$. The cost is approximately $4b_{L-1}Ms = (N-4s)M < MN$ flops per parity.

For all the lower levels the matrix vector product in Eq. \eqref{eq:wTf} is computed more efficiently through Step 2 in Alg. \ref{alg:fastercomplete}. This step represents a novel contribution of the current work, and we shown is some details in Appendix \ref{sec:A1} how Eq. \eqref{eq:wTf} can simply be computed as 
\begin{equation}
    \boldsymbol{w}(\gamma, b, q) = \sum_{l\in \{0, 1\}} \boldsymbol{B}{(l)}\boldsymbol{w}(\gamma+1, 2b+q+1, l), \label{eq:wtk}
\end{equation}
where $\boldsymbol{B}(l)=(b^{(l)}_{kn})_{k,n=0}^{M-1}$ are the lower triangular matrices that map the basis function $T_k$ used on the half-interval $I(l)$ of level $\gamma$ (see illustration for matrix $A(\gamma, b, 0)$ in Fig. \ref{fig:2levelsFMM}) to the full-interval Chebyshev polynomials used for the matrices directly below  on level $\gamma+1$. That is
\begin{equation}
    T_k((Y+2l-1)/2) = \sum_{j=0}^k b^{(l)}_{kj} T_j(Y), \quad Y \in I, \, l \in \{0, 1\}, \label{eq:Tmap}
\end{equation}
where $(Y+2l-1)/2 \in I(l)=[l-1,l]$. The two lower triangular matrices $\boldsymbol{B}(l) \in \mathbb{R}^{M \times M}$ can be computed exactly using the binomial theorem and simple basis changes, see Appendix \ref{sec:A3}. Equation \eqref{eq:Tmap} is the modal version of Eq. (4.9) used in \cite{alpert91}. However, the corresponding matrices (i.e, $(u_r(t_i/2))_{r,i=0}^{M-1}$ in Eq. (4.9) of \cite{alpert91}) required by the Lagrange interpolating polynomials are both dense.
  
Equation \eqref{eq:wtk} consists of the sum of two matrix vector products, where the matrices $\boldsymbol{B}(l) \in \mathbb{R}^{M \times M}$ are lower triangular. The cost would thus normally be $2M^2$ flops for each parity and nearly twice as fast as the dense approach used by Alpert and Rokhlin \cite{alpert91}. However, we can also take advantage of the fact that the two matrices $\boldsymbol{B}{(0)}$ and $\boldsymbol{B}{(1)}$ are very similar (see Eq. \eqref{eq:appT} and App. \ref{sec:A3}, $\boldsymbol{B}(0)$ and $\boldsymbol{B}(1)$ differ only in the sign of the odd diagonals), and the cost may then be reduced even further to one single lower triangular matrix vector product ($M^2$ flops) plus a small additional cost of $2M$ for preparing two vectors of length $M$. The flop count is thus nearly 4 times less than a dense approach. A faster algorithm is described in detail in Alg. \ref{alg:Tw} in Appendix \ref{sec:A2}. 

Step (3) of Alg.~\ref{alg:fastercomplete} is simply the application of the precomputed coefficient matrices $\boldsymbol{\hat{A}}(\gamma, b, r)$  to the computed work vector $\boldsymbol{w}(\gamma, b, q)$ on each submatrix. The cost is here the same $4N_c M^2<6NM^2/s$ as computed in Eq. \eqref{eq:totalnumflops2}.

Step (4) of Alg.~\ref{alg:fastercomplete} is the second matrix vector product from Alg.~\ref{alg:matvecblock} ($\boldsymbol{z} = \boldsymbol{{T}} \boldsymbol{c}(p)$) that needs to be handled more efficiently in order to obtain a $\mathcal{O}(N)$ method. Since this is basically the transpose of Step (2) we can also here make use of the same matrices $\boldsymbol{B}(l)$. However, the computed result will now have to be transported in the reverse direction, from the lowest level 0 and up to the highest level $L-1$. 
Since $\boldsymbol{B}(l)$ is triangular the cost for the two matrix vector products in the innermost loop of step (4) will be $2M^2$ flops per parity. This can be further improved to $\approx M^2$ flops by taking advantage of the structure of $\boldsymbol{B}(0)$ and $\boldsymbol{B}(1)$, similar to Step (2). 

Having transported the intermediate array $\boldsymbol{c}$ up to the highest level, $L-1$, a tensor contraction (Step (5) in Alg. \ref{alg:fastercomplete}) of the same cost as Step (1) concludes the FMM part of the method
\begin{equation}
\boldsymbol{z}^{\sigma}(L-1) = \boldsymbol{c}(L-1) {\boldsymbol{T}^T(\sigma, L-1)}. \label{eq:step5}
\end{equation}

The complete and faster $\mathcal{O}(N)$ Legendre to Chebyshev (L2C) algorithm presented in Alg. \ref{alg:fastercomplete} is basically 5 different matrix vector products (or tensor contractions) and simple enough that it can be implemented in less than 100 lines of high-level Python code. The algorithm for the reverse direction (C2L) simply requires that the diagonal matrix $\boldsymbol{C}$ is applied to the input vector first (see Eq. \eqref{eq:directmv}), and otherwise the only difference is that $\mathcal{L}(x,y)$ (see Eq. \eqref{eq:Lxymod} and Eq. (2.21) of \cite{alpert91})
is used instead of $\mathcal{A}(x,y)$ for the Chebyshev approximations in Eqs. \eqref{eq:akl2} and \eqref{eq:aklmore}.

\begin{algorithm}
\caption{Complete faster algorithm for the Legendre to Chebyshev transform (L2C) $\boldsymbol{z} = \boldsymbol{C}^{-1} \boldsymbol{A} \boldsymbol{f}$.}    \label{alg:fastercomplete}
\begin{algorithmic}
    \Function{Leg2Cheb}{$L,\boldsymbol{f}, \boldsymbol{T}, \boldsymbol{\hat{A}}, \boldsymbol{B}, \boldsymbol{C}$}
    \Input{$L$}{$\mathbb{N}$}{Number of levels}
    \Input{$\boldsymbol{{f}}$}{$\mathbb{R}^{N}$}{Legendre coefficients}
    \Input{$\boldsymbol{{T}}$}{$\mathbb{R}^{2 \times s \times M}$}{Vandermonde matrices, $s=h_{L-1}$}
    \Input{$\boldsymbol{\hat{A}}$}{$N_c$ sequence of $\mathbb{R}^{M \times M}$}{Chebyshev matrices of coefficients}
    \Input{$\boldsymbol{B}$}{$\mathbb{R}^{2 \times M \times M}$}{Lower triangular matrices}
    \Input{$\boldsymbol{C}$}{matrix $\in \mathbb{R}^{N \times N}$}{diagonal matrix, see Eq. \eqref{eq:AfCf}}
    \Output{$\boldsymbol{{z}}$}{$\mathbb{R}^{N}$}{Chebyshev coefficients}
    \For{$\sigma \gets 0, 1$}
        \For{$\gamma \gets 0, L-1$} \Comment{Initialize work arrays}
            \State $\boldsymbol{c}(\gamma) \gets 0$ \Comment{Work array  $\mathbb{R}^{b_{\gamma} \times 2 \times M}$}
            \State $\boldsymbol{w}(\gamma) \gets 0$ \Comment{Work array $\mathbb{R}^{b_{\gamma} \times 2 \times M}$}
        \EndFor
        \State $\boldsymbol{w}(L-1) \gets \boldsymbol{f}^{\sigma}(L-1) \boldsymbol{T}(\sigma)$ \Comment{See Eq. \eqref{eq:wMM}. Step (1)}
        \For{$\gamma \gets L-1, 1$}\Comment{See Eq.~\eqref{eq:wtk}. Step (2)}
            \For{$b \gets 1, b_{\gamma}-1$}
                \State $b_0, q_0 \gets \textbf{divmod}(b-1, 2)$
                \State $\boldsymbol{w}(\gamma-1, b_0, q_0) = \boldsymbol{B}(0)\boldsymbol{w}(\gamma, b, 0)+ \boldsymbol{B}(1)\boldsymbol{w}(\gamma, b, 1)$  
            \EndFor
        \EndFor
        \For{$\gamma \gets 0, L-1$} \Comment{Step (3)}
            \For{$b \gets 0, b_{\gamma}-1$}
                \For{$q \gets 0, 1$}
                    \For{$p \gets 0, q$}
                        \State $r \gets p + q(q+1)/2$ \Comment{See Eq. \eqref{eq:cfrompq}}
                        \State $\boldsymbol{{c}}(\gamma, b, p) \gets  \boldsymbol{{c}}(\gamma, b, p) + \boldsymbol{\hat{A}}(\gamma, b, r )\boldsymbol{w}(\gamma, b, q)$
                    \EndFor
                \EndFor
            \EndFor
        \EndFor
        %
        \For{$\gamma \gets 0, L-2$} \Comment{Step (4)}
            \For{$b \gets 0, b_{\gamma+1}-2$}
                \State $b_0, q_0 \gets \textbf{divmod}(b, 2)$
                \For{$p \gets 0, 1$}
                    \State $\boldsymbol{{c}}(\gamma+1, b, p) \gets  \boldsymbol{{c}}(\gamma+1, b, p) + \boldsymbol{B}^T(p) \boldsymbol{c}(\gamma, b_0, q_0)$
                \EndFor
            \EndFor
        \EndFor
        \State $\boldsymbol{z}^{\sigma}(L-1) \gets \boldsymbol{c}(L-1) {\boldsymbol{T}^T(\sigma)}$ \Comment{See Eq. \eqref{eq:step5}. Step (5)}
    \EndFor
    \State Apply direct method to remaining part \Comment{Step (6)}
    \State $\boldsymbol{z} \gets \boldsymbol{C}^{-1} \boldsymbol{z}$ \Comment{Step (7)} \\
    \Return $\boldsymbol{z}$
    \EndFunction
\end{algorithmic}
\end{algorithm}

\subsection{Total cost of faster method}
\label{sec:total_cost_faster}
\begin{lemma}
Retaining only the most deeply nested loops, the total flop count of the faster method described in Alg. \ref{alg:fastercomplete} can be estimated as
\begin{equation}
    \left( \mu s + 4M +  \frac{(6+2) M^2}{s}\right)N. \label{eq:totalnumflopsnew}
\end{equation}

\begin{proof}
    The algorithm is divided into 7 steps, as shown in Alg. \ref{alg:fastercomplete}.
    There is no difference in cost from \eqref{eq:totalnumflops2} for step (6): $\mu s N$ or step (3): $6NM^2/s$. For steps (1) and (5) the tensor contractions each cost less than $MN$ flops (see Eq. \eqref{eq:wMM}) for each parity and thus a total of $4MN$. The two matrix vector products in the innermost loop of both steps (2) and (4) cost for each parity $2M^2$ plus a small and neglected cost for setting up vectors of length $M$. These products are applied for all but one block on all levels, which counts to  $\sum_{\gamma=1}^{L-1}(b_{\gamma}-1)=2(2^L-L-1)<\frac{N}{2s}$ blocks (see Eqs.~\eqref{eq:total_blocks_and_mats} and \eqref{eq:Nb}). Steps (2) and (4) thus costs $\approx 4M^2 \cdot {N}/{2s}=2M^2N/s$ flops. Step (7) is simply $N$ multiplications and thus neglected along with similar small operations.
\end{proof}
\end{lemma}

It is noteworthy that steps (2) and (4) now together cost approximately $2M^2N/s$ flops and much less than step (3) (i.e., $6M^2N/s$). For the nodal Lagrange approach steps (2) and (4) require approximately $8M^2N/s$ flops. 

Assuming that computational speed is simply proportional to the number of flops, we can compute the optimal size of the smallest submatrices from Eq. \eqref{eq:totalnumflopsnew}. Setting $\mu=3$ and $M=18$, we find that the optimal $s$ is $29$, and for a power of 2 input this should be rounded up to $s=32$, which is the number used by Alpert et al. \cite{alpert91}. For these parameters the modal L2C requires approximately $250N$ flops for execution, which is nearly $20~\%$ less than the nodal L2C ($310N$). We note that a DCT type 2 costs $\approx 2N\log_2 N$  flops (see \cite{Shao2008}) and the modal FMM is thus using $125 / \log_2 N$ times more flops than a DCT.

Finally, we note that the memory requirement for the method described in Alg. \ref{alg:fastercomplete} with $s=32$ is merely $ \approx 17N$ numbers of the chosen double precision, if we choose not to precompute $\boldsymbol{A}$ for the direct part. This is the requirement for storing $N_c$ matrices $\boldsymbol{\hat{A}}$, each of shape $\mathbb{R}^{M \times M}$, and two vectors of shape $\mathbb{R}^N$ (a work array and $\boldsymbol{\lambda}$). Everything else is negligible in terms of memory (assuming $N\gg s$). If we choose to precompute entries of $\boldsymbol{A}$ for the direct part, there is an additional and significant memory requirement of $1.5sN$ numbers.

\subsection{A faster initialization}
\label{sec:fastinit}
The major cost for the initialization of the fast multipole method is to run over all submatrices, evaluate (see Eq. \eqref{eq:aklmore})
\begin{equation} 
\mathcal{A}(x_m, y_n) = \Lambda \left(\frac{y_n-x_m}{2}\right) \Lambda \left(\frac{y_n+x_m}{2} \right), \quad 0 \le m, n < M,
\end{equation}
and apply a two-dimensional DCT of type 2 to each matrix $(\mathcal{A}(x_m, y_n))_{m,n=0}^{M-1}$. We note that the matrix $(\mathcal{A}(x_m, y_n))_{m,n=0}^{M-1}$ is the Hadamard product of two matrices, that for simplicity will be referred to in this section as $\Lambda^- = (\Lambda((y_n-x_m)/2))_{m,n=0}^{M-1}$ and $\Lambda^+ = (\Lambda((y_n+x_m)/2))_{m,n=0}^{M-1}$. 

For the initialization it is very important that $\Lambda(z)$ is computed efficiently. Alpert et al.~\cite{alpert91} describe how asymptotically $\Lambda(z) \sim \frac{1}{\sqrt{z}}$ as $z \longrightarrow \infty$ such that $\Lambda(z)\sqrt{z}$ can be well approximated by a polynomial series in $1/z$ of length 6. An even better approach along the same lines is suggested by Bogaert et al. \cite{bogaert12} and  Slevinsky \cite{slevinsky2019}, which is to adapt a Taylor series in $1/z$ for the function $\tau(z) = \sqrt{z}\frac{\Gamma(z+\sfrac{1}{4})}{\Gamma(z+\sfrac{3}{4})}$ and then compute $\Lambda(z) = \tau(z+\sfrac{1}{4})/\sqrt{z+\sfrac{1}{4}}$. This latter approach is better because the Taylor expansion is more quickly decaying, with only even terms remaining:
\begin{equation}
    \tau(z) \approx 1 - \frac{1}{2^6 z^2} + \frac{21}{2^{13}z^4} - \frac{671}{2^{19}z^6}+ \frac{180,323}{2^{27}z^8} +\mathcal{O}\left(\frac{1}{z^{10}}\right). \label{eq:tau}
\end{equation}
Another advantage is that all but the last of these rational coefficients can be represented exactly in double precision since the denominator of the coefficient is small enough and a power of 2. Using Eq. \eqref{eq:tau} we find that $\Lambda(z)$ can be computed with relative error to within machine precision (i.e., $\approx 2.22\times 10^{-16}$) for all $z > 19.88$. However, as $z\longrightarrow \infty$ more terms on the right hand side can be discarded and it is easy to verify numerically that with only $4, 3$ or 2 terms, the expression is still accurate to machine precision for $z > 40, 134$ or $1844$, respectively.\footnote{Here rounding up to nearest integer.} Hence, the cost for evaluating $\mathcal{A}(x_m, y_n)$ will be smaller for large $N$. 
The smallest relevant $(x, y)$ for the submatrices are found on the highest level $L-1$ closest to the main diagonal where $(y-x)/2 \ge s$. Hence, with $s=64$, we only need the first 4 terms of Eq. \eqref{eq:tau} on the highest level $L-1$. For levels $L-2, L-3, L-4$ and $L-5$ it is sufficient to use 3 terms and for all levels $<L-5$ merely 2 terms of Eq.~\eqref{eq:tau} are required. We can estimate that $\Lambda(z)$ costs $20$ flops for the highest level (including $\approx 10$ for a square root), and then remove 3 flops per removed term of Eq.~\eqref{eq:tau} for the higher levels. 

There are still significant enhancements that can be made for more efficient evaluation of the matrix $(\mathcal{A}(x_m, y_m))$, by taking advantage of certain symmetries. First of all, the matrix $\Lambda^-$ is persymmetric, whereas $\Lambda^+$ is symmetric. Hence, symmetry dictates that only $(M^2+M)/2$ items of each matrix need to be evaluated, which nearly halves the flop count of evaluation. Second, since $\Lambda^-$ is persymmetric and only depends on the distance to the main diagonal of the matrix $\boldsymbol{A}$ (see Eq. \eqref{eq:hik}), it is sufficient to evaluate $\Lambda^-$ for the first two submatrices of the first block of any given level (i.e., for the submatrices $\boldsymbol{A}(\gamma, 0, 0)$ and $\boldsymbol{A}(\gamma, 0, 1)$, see Eq. \eqref{eq:Aij}). All remaining submatrices on that level will then be able to reuse one of these two already evaluated matrices. Asymptotically this eliminates the cost for evaluating $\Lambda^-$, as all it takes now is $L(M^2+M)$ evaluations, with $L \sim \log_2 N$, see Eq. \eqref{eq:numlevels}. Since computing $\Lambda^-$ is half the cost of evaluating $(A(x_m, y_n))$ the flop count is thus reduced by nearly another factor 2. With these optimizations the flop count for the evaluation drops nearly by a factor 4 from the na\"ive evaluation of $\mathcal{A}(x_m, y_n)$ at all points. Assuming $\Lambda(z)$ costs 20 flops, then, on average, the cost will be $10M^2$ ($40M^2$ for na\"ive, divided by 4) flops for each submatrix.

For the methods described in Secs. \ref{sec:FMM} and \ref{sec:faster} we also need to take the DCT of all the matrices $(\mathcal{A}(x_m, y_n))$. The cost of this DCT will depend on the implementation, but it should be close to $2M \log_2 M$, see \cite{Shao2008}. We have implemented the DCT of size $M=18$ using a mixed radix 2 and 3, since $18 = 2 \cdot 3^2$, at a cost close to $2M \log_2 M$ flops. Hence, the 2D DCT on a submatrix, where the 1D DCT is applied to $2M$ vectors of length $M$, can be computed in approximately $4M^2\log_2 M$ arithmetic operations.

For the direct part of the transform (step 6) we need to precompute the vector $\boldsymbol{\lambda}=\{\Lambda(i)\}_{i=0}^{N-1}$. A direct application of Eq. \eqref{eq:tau} at 20 flops per evaluation will here lead to an initialization cost of $\approx 20 N$ flops. Some of these could be stored in a lookuptable, but for integer arguments to $\Lambda(z)$ a more elegant solution is to use $\Lambda(0)=\Gamma(\sfrac{1}{2})=\sqrt{\pi}$ and then the recursion
\begin{equation}
    \Lambda(i) = \Lambda(i-1)\frac{i-\sfrac{1}{2}}{i}, \quad \forall \, i=1, 2, \ldots, N-1, \label{eq:Lambdaint}
\end{equation}
which is stable since $(i-\sfrac{1}{2})/i$ remains close to, and is always smaller than, 1. The absolute error of computing Eq. \eqref{eq:Lambdaint} is found to be less than machine precision for all $i < 2^{23}$. However, the error of each iteration is found to be 2 ulps (unit in the last place) or less, and these may accumulate. In order to achieve a relative error for all $i<2^{23}$ that is less than $ 10^{-15}$ and an absolute error less than 5 ulps, we have found it necessary to compute every $8$'th $\Lambda(i)$ from Eq. \eqref{eq:tau} and the remaining 7 out of 8 from \eqref{eq:Lambdaint}. The cost for computing $\boldsymbol{\lambda}$ will then be approximately $\tfrac{7}{8}2N + \tfrac{1}{8}20N \approx 4 N$ flops. If all items of the matrix $\boldsymbol{A}$ that are used in the direct part are precomputed, then the initialization cost will be $4N + 1.5sN$ flops, since in addition to computing $\boldsymbol{\lambda}$, we would also need to multiply two $\Lambda(i)$'s together for each of the $1.5sN$ nonzero items. However, since it has been found faster to not use precomputation, and since the memory requirement is many times larger and the initialization takes much longer, we will only consider the non-precomputed version of the algorithm in what follows.

In summary, the cost for the Legendre to Chebyshev (L2C) initialization may be estimated as a small fixed cost plus $N_c (4M^2 \log_2 M + 10 M^2)+4N$ flops asymptotically. Inserting for $N_c$ we obtain $\approx  N(6M^2 \log_2 M+15M^2)/s+4N$ flops.  For the nodal Lagrange method the DCT, and thus the term $6M^2 \log_2 M$, can be dropped.

Finally, we note that the transform in the reverse direction, Chebyshev to Legendre (C2L), can be computed without any further evaluations of the $\Lambda(z)$ function. This follows since the transform function for the reverse direction, $\mathcal{L}(x, y)$, can be written such as to reuse $\Lambda^-$ and $\Lambda^+$ for its nonzero items:
\begin{equation}
    \mathcal{L}(x, y) = \frac{2(x+\sfrac{1}{2})y}{(x+y)(x+y+1)(x-y+1)}\frac{\Lambda \left( (y-x)/2\right)}{\Lambda \left( (x+y)/2\right)}.
    \label{eq:Lxymod}
\end{equation}
Equation \eqref{eq:Lxymod} is a simple reformulation of Eq. (2.21) of \cite{alpert91}, using Eq. \eqref{eq:Lambdaint} and $\Lambda(i-\sfrac{1}{2})=(i \Lambda(i))^{-1}$. 


\section{Results}
\label{sec:results}
All numerical results are computed on a MacBook M3 Pro (4.05 GHz 12 Cores) and/or the SAGA supercomputer (Intel Xeon-Gold 6138 2.0 GHz) provided by 
Sigma2 - the National Infrastructure for High Performance Computing and Data Storage in Norway. The algorithms have been implemented in both C and Python, and the code is available in the public repository \url{https://github.com/mikaem/SISC-Legendre-to-Chebyshev.git}. The Python code is not very efficient and mainly intended to illustrate the simplicity of the algorithm. The C-code is more optimized for speed, with the most heavy duty operations outsourced to a BLAS vendor like OpenBLAS \cite{openblas} (SAGA) or Accelerate \cite{accelerate} (Mac). The triangular matrix-vector products in steps (2) and (4) and the direct step (6) are the only parts that are implemented in detail by hand. All computations use constant $M=18$ and $s=32$ since this has been found optimal in terms of both flops and execution speed. These are the same parameters as used by Alpert and Rokhlin \cite{alpert91}.\footnote{There is a factor 2 difference between our definition of $s$ and the parameter $s$ used by \cite{alpert91}. This is because we need the submatrix-sizes to be even and thus define the smallest submatrices to have size $2s$ instead of $s$.} The direct part (step (6)) of Alg. \ref{alg:fastercomplete} is performed without precomputing matrix entries. All computations are performed at first using a single core. We consider multithreading in the last paragraph of this section.

\subsection{Accuracy}
To evaluate the numerical accuracy we use the relative maximum norm
\begin{equation}
    E_{\infty} = \frac{||\boldsymbol{z} - \boldsymbol{z}^*||_{\infty}}{||\boldsymbol{z}^*||_{\infty}}, \label{eq:Einf}
\end{equation}
where $\boldsymbol{z}$ is the computed result and $\boldsymbol{z}^*$ is the 'exact' result computed with a direct multiprecision C++ solver, accurate to approximately 32 digits. We will also evaluate two consecutive forward and backward transforms, in which case the exact result will be the input data.

We use random coefficients that decay as $N^{-r}$, where $r \ge 0$, see, e.g., \cite{shen19,Hale2014,Townsend18}. The random coefficients are drawn uniformly from the interval $[0, 1]$ using the standard C function $\textbf{rand}$ and for repeatability we use the constant seed $\textbf{srand(1)}$. At first we consider only the accuracy of the Legendre-to-Chebyshev (L2C) and Chebyshev-to-Legendre (C2L) methods for uniform coefficients with $r=0$. Since the multiprecision implementation uses a direct method we limit the test for coefficient vectors of size up to $N= 2^{16}$. Table \ref{tab:error1} shows that the algorithm is accurate to all but one digit for both directions and there is very little error growth with size. Also, there is hardly any difference in the $E_{\infty}$ norm between the modal approach described in this paper and the original nodal approach of \cite{alpert91}. The apparent coincidence that most $E_{\infty}$ numbers are the same for both nodal and modal approaches in Table \ref{tab:error1} is attributed to the fact that the error $||\boldsymbol{z} - \boldsymbol{z}^*||_{\infty}$ is often merely a small integer or half-integer times one ulp of $||\boldsymbol{z}^*||_{\infty}$ (the value of the ulp depends on $||\boldsymbol{z}^*||_{\infty}$). The error in $||\boldsymbol{z} - \boldsymbol{z}^*||_{\infty}$ in ulps of $||\boldsymbol{z}^*||_{\infty}$ is given in parenthesis in Table \ref{tab:error1}. 
\begin{table}[]
    \centering \small
    \caption{Maximum error $E_{\infty}$ of Legendre-to-Chebyshev (L2C) and Chebyshev-to-Legendre (C2L) transforms for the current modal (m) FMM and the nodal (n) FMM of \cite{alpert91}. The input array is drawn from random numbers in the range $[0, 1]$ and the error is computed with \eqref{eq:Einf} using a direct multiprecision solver for $\boldsymbol{z}^*$. The number in parenthesis is the error in $||\boldsymbol{z} - \boldsymbol{z}^*||_{\infty}$ as the number of ulps of $||\boldsymbol{z}^*||_{\infty}$.}
\begin{tabular}{ccccc}
 N & L2C (m) & L2C (n) & C2L (m) & C2L (n) \\ 
\hline 
256 & 8.88e-16 (4.0) & 8.88e-16 (4.0) & 7.44e-15 (4.2) & 7.44e-15 (4.2) \\ 
512 & 1.11e-15 (2.5) & 8.88e-16 (2.0) & 1.10e-14 (3.1) & 1.10e-14 (3.1) \\ 
1024 & 1.11e-15 (2.5) & 1.33e-15 (3.0) & 2.16e-14 (6.1) & 2.16e-14 (6.1) \\ 
2048 & 1.11e-15 (2.5) & 8.88e-16 (2.0) & 3.91e-14 (11.0) & 3.91e-14 (11.0) \\ 
4096 & 2.44e-15 (5.5) & 1.55e-15 (3.5) & 5.68e-14 (8.0) & 5.68e-14 (8.0) \\ 
8192 & 1.78e-15 (4.0) & 2.22e-15 (5.0) & 9.59e-14 (13.5) & 9.59e-14 (13.5) \\ 
16384 & 2.44e-15 (5.5) & 2.44e-15 (5.5) & 1.39e-13 (9.8) & 1.39e-13 (9.8) \\ 
32768 & 2.44e-15 (5.5) & 2.22e-15 (5.0) & 1.99e-13 (14.0) & 2.06e-13 (14.5) \\ 
\hline 
\end{tabular}

\label{tab:error1}
\end{table}

\begin{figure}
    \centering
    \includegraphics[width=0.6\textwidth]{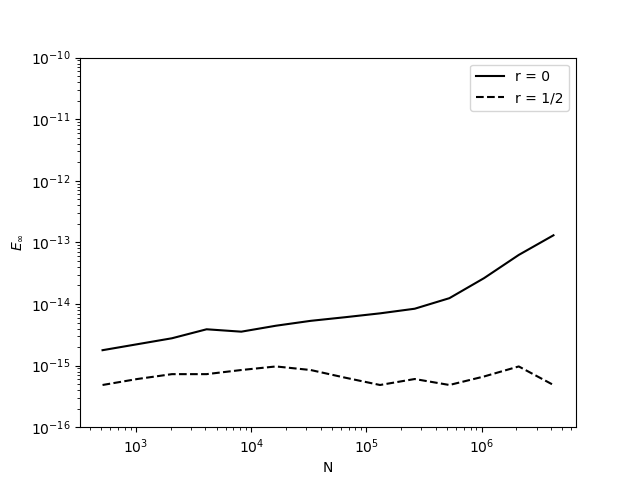}
    \caption{The normalized maximum error for a consecutive L2C and C2L transforms. The black solid curve is using uniform random input data (r=0) and the dashed curve is using decaying, random  Legendre coefficients scaled with r=1/2.}
    \label{fig:error}
\end{figure}

In Fig. \ref{fig:error} we show the error of computing one L2C followed by a C2L, which should return the input array. For this test we use both uniform random coefficients (r=0) and coefficients decaying with r=1/2. The latter is more realistic since Legendre (or Chebyshev) coefficients usually decay in applications where the solution is smooth. We use only the modal FMM since the difference in accuracy from the nodal is minimal. From Fig. \ref{fig:error} we see that the error for the decaying input array is close to machine precision (one ulp for this plot equals $1.11 \times 10^{-16}$) for $N$ all the way up to 10 million, whereas the uniform input leads to slightly decreased accuracy that becomes notable for $N>100,000$.  Similar bounded $\mathcal{O}(1)$ accuracy was also obtained by the Toeplitz-Hadamard (TH) method of Townsend et al. \cite{Townsend18} for decay rates r $\ge$ 1/2, whereas they observed $\mathcal{O}(N^{0.5}\log N)$ error growth for r=0. However, they only presented results for $N<1000$. Results for $N$ up to $10^7$ were presented for the divide and conquer model by Olver et al. \cite{Olver2020}, showing a $\mathcal{O}(1)$ error for r=1. However, their simulation with $N=8388608$ required $10^{11}$ bytes of memory, see Table 4.3 of \cite{Olver2020}. By comparison, the current method requires approximately $10^9$ bytes ($17N$ double precision numbers) for the same problem size. Finally, we note that the error growth was considerably larger both in Hale et al. \cite{Hale2014} and Shen et al. \cite{shen19}, whereas Keiner \cite{keiner09}, who used the original FMM method, could also show $\mathcal{O}(1)$ error growth for r=0 and $N \le 10^5$ (see Table 6.2 of \cite{keiner09}). 

\subsection{Computational speed, single-core}
\begin{figure}[t]
    \centering
    \begin{subfigure}[t]{0.49\textwidth}
        \includegraphics[width=\textwidth]{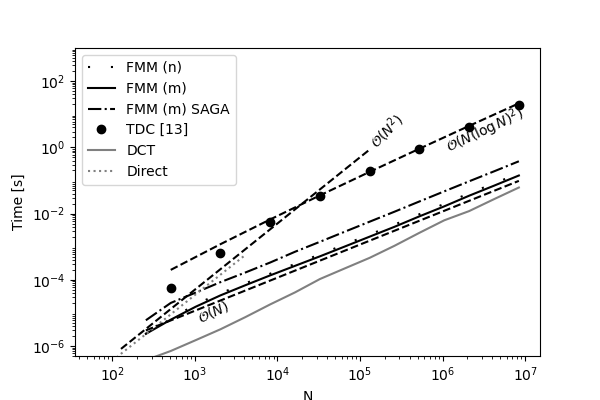}
        \caption{Execution times.}
    \end{subfigure}
    \begin{subfigure}[t]{0.49\textwidth}
        \includegraphics[width=\textwidth]{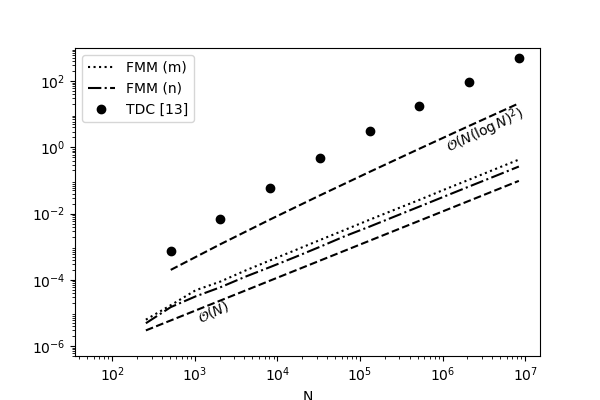}
        \caption{Planning times.}
    \end{subfigure}
    \caption{In (a) the black solid and dash-dotted curves show the execution times in seconds for a Legendre to Chebyshev transform of the current FMM (m) on the MacBook and SAGA computers, respectively. The loosely dotted curve shows the execution time (Mac) for the nodal FMM (n). The large black dots show the execution times of the triangular-banded divide and conquer (TDC) method from Table 4.3 of \cite{Olver2020}. The execution time (Mac) for a type 2 DCT (FFTW) is shown in gray. In (b) the planning time (Mac) for the FMM (m) is shown as a black dotted line, for FMM (n) in a dash-dotted line and for TDC \cite{Olver2020} in large dots. Three dashed lines are shown with increasing slope to illustrate the asymptotic behaviour for $\mathcal{O}(N), \mathcal{O}(N (\log N)^2)$ and $ \mathcal{O}(N^2)$. The asymptotic curves are identical in (a) and (b).}
    \label{fig:speed}
\end{figure}

The speed (wall clock time) of the current modal (single-core) L2C\footnote{The speed of the reverse direction C2L is basically identical, so only one direction is shown.} transform is shown in Figure \ref{fig:speed} (a) with a black solid curve for the MacBook and a dash-dotted curve for the SAGA supercomputer. The large black dots show the data reported in Table 4.3 of \cite{Olver2020} for the triangular-banded divide and conquer (TDC) method. The gray solid curve shows the execution time of a single-core DCT type 2 computed with precompiled (conda-forge, using --enable-sse2 and --enable-avx) FFTW \cite{fftw} on the MacBook, using the $\text{FFTW\_MEASURE}$ planning flag. All measurements are presented as the fastest executions obtained for $2^{26-\log_2 N}$ repetitions.  The $\mathcal{O}(N)$ scaling is obvious for the current FMM method on both the SAGA computer and on the Mac. It is noteworthy that a transform the size of $N=10^7$ can be obtained in approximately 0.2 seconds with the current FMM, which is approximately three times slower than the type 2 DCT computed with FFTW. 
For large $N$ the speed of FMM can bee seen to be nearly 2 orders of magnitude faster than the execution times reported for the TDC method \cite{Olver2020}. A large part of this factor may be due to implementation, but the results of TDC are also reported to be faster than the comparative Toeplitz-Hankel method described in \cite{Townsend18} (see Fig 5 (right)), which again are reported to be faster than the method of Hale et al. \cite{Hale2014}. The speed of the nodal FMM using Lagrange polynomials is shown as a loosely dotted curve, and it falls very close to the results of the modal FMM, being on average $\approx 15-20 \%$ slower. The slower execution is only due to the transport of intermediate results through the novel steps (2) and (4) in Alg. \ref{alg:fastercomplete}. Otherwise, the execution makes use of exactly the same code, only with different matrices and coefficients.
A direct comparison with execution times of Shen et al. \cite{shen19} is not possible, because they only report the flop count of their execution. However, comparing to Fig 3 in \cite{shen19}, the current minimal flop count of $\approx 250N$ is at least a factor 2 lower.


In Fig. \ref{fig:speed} (b) we show the time spent on planning the FMM transforms, both for the current method (black dotted line) and for the Lagrange FMM method (dash-dotted line) using the optimizations described in Sec \ref{sec:fastinit} without precomputation of matrix entries. The TDC from Table 4.3 of \cite{Olver2020} is shown as large black dots. We observe the claimed $\mathcal{O}(N)$ scaling, and a  modal FMM planning that takes approximately three times longer than execution. The optimized planning stage for the nodal FMM is very nearly as fast as the execution. We note that the nodal FMM \cite{alpert91} is faster to plan since it does not compute the DCT on its data. Otherwise the planning is very similar. Also note that using a larger value for $s$, like $s=64$, the planning will become faster, whereas the execution becomes slower. This is because there are then fewer submatrices to initialize, whereas the direct part of the execution becomes more expensive.

\subsection{Computational speed, multi-core}
\begin{figure}
    \centering
    \includegraphics[width=0.6\textwidth]{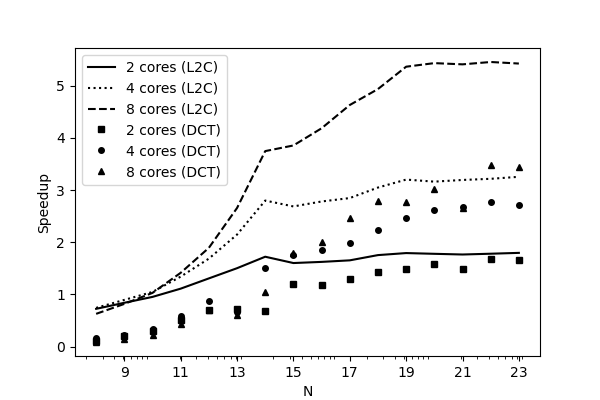}
    \caption{Speedup measured for the L2C method is shown in black solid, dotted and dashed lines for 2, 4 and 8 CPU cores, respectively. For FFTW's DCT type II the corresponding measurements are shown as squares, circles and triangles. }
    \label{fig:threads}
\end{figure}

Algorithm \ref{alg:fastercomplete} lends itself easily to shared memory parallelization since it can make use of multithreaded BLAS routines in order to compute its numerous matrix-matrix and matrix-vector products. This is the obvious and only approach for steps 1 and 5. For steps 2, 3 and 4, on the other hand, we have found it more efficient to use single-core matrix operations and instead distribute threads over blocks. That is, we place the OpenMP directive '\lstinline{#pragma omp parallel for}' just prior to the three for-loops running over blocks in steps 2, 3 and 4 in Alg. \ref{alg:fastercomplete}. The direct step 6 has also been implemented using a loop over blocks that is parallelized with OpenMP. We compute the multi-core results on the SAGA supercomputer, using up to 8 CPU cores. We utilize precompiled cluster modules for both OpenBLAS and FFTW, both compiled and linked with support for OpenMP. We compare the parallel performance for the execution speed of L2C with a DCT type II from FFTW, because of the DCT's relevance to complete Legendre transforms (see, e.g., \cite{Hale2014}). Figure \ref{fig:threads} shows the speedup gained by both methods for a range of input array sizes, using 2, 4 or 8 CPU cores. The results shown are the fastest executions obtained from 100 repetitions for each size.  We note that the DCT requires quite large arrays ($>2^{14}=16384$) in order to obtain a speedup, whereas the L2C shows speedup for arrays $>2^{10}=1024$. For all array sizes the speedup of L2C for a given number of cores is greater than for the DCT. We note that perfect scaling would correspond to a speedup of $n$ for $n$ CPUs and the performance depends heavily on both the computer's hardware and the compilation of the software. 

\section{Conclusions}
In this paper we have described a modal version of the originally nodal Fast Multipole Method (FMM) of Alpert and Rokhlin \cite{alpert91}, used for transforming between Legendre and Chebyshev coefficients. The modal approach allows us to spread intermediate results faster through different levels of hierarchical matrices, allowing us to speed up the original method with approximately $20~\%$. For a power of 2 double precision input array, we can perform the modal transform in approximately $250N$ flops, compared to $310N$ for the nodal approach. However, the modal approach is slower to plan, due to an additional discrete cosine transform required by the initialization. The accuracy is the same for both nodal and modal versions, and we can transform millions of slightly decaying coefficients at nearly machine precision.

We have described in detail an algorithm that is easy to implement, relying mostly on a series of dense matrix-vector products that can be executed efficiently through BLAS. A single-core C implementation is shown to be able to transform $N=10^6$ double precision coefficients in approximately $20$ milliseconds on a MacBook M3 Pro laptop. This is approximately three times slower than a well-planned (single-core) discrete cosine transform of type 2 from FFTW \cite{fftw} on the same machine.  For large $N$ the planning of the transform, which needs to be done just once for given $N$, is approximately 3 times slower than the execution, and the memory used for a double precision transform is no more than $17 N$ numbers. We have also described an efficient OpenMP implementation as well as optimizations for the nodal FMM, which lead to at least 4 times faster initialization than for the original nodal method. 

The method described in this paper is easily adapted to similar transforms between any bases of Gegenbauer polynomials, and not only Legendre-Chebyshev. For general Jacobi polynomials, that are not necessarily odd or even functions, the method can be extended by using a full upper triangular connection matrix. 

\section{Code availability}
The code used to implement the method described in this paper is open source and available from the repository \url{https://github.com/mikaem/SISC-Legendre-to-Chebyshev.git}. Here both the original nodal FMM and the new modal FMM are implemented. The method is already in use in the spectral Galerkin framework Shenfun \cite{mortensen_joss}. 

\section{Acknowledgements}
Some computations were performed on resources provided by 
Sigma2 - the National Infrastructure for High Performance Computing and Data Storage in Norway. 

\bibliography{bib.bib}

\begin{appendices}

\section{Appendix}


\renewcommand{\theequation}{A-\arabic{equation}}%
\vskip 0.1in

\subsection{Derivation of Eq. \eqref{eq:wtk}}
\label{sec:A1}

In order to derive Eq. \eqref{eq:wtk} we will be using the same variables and parameters as defined in Sec. \ref{sec:faster}. We start by assuming that the work array $\boldsymbol{w}(\gamma+1) \in \mathbb{R}^{b_{\gamma+1} \times 2 \times M}$, defined prior to Eq. \eqref{eq:wTf}, is known, and the objective is to efficiently compute $\boldsymbol{w}(\gamma)$. To this end we first write Eq. \eqref{eq:wTf} in index form for each parity as
\begin{equation}
    w_n(\gamma, b, q) = \sum_{m=0}^{h_{\gamma}-1} T_n (Y^{\sigma,\gamma}_m) {f}_m^{\sigma}(\gamma, b, q), \quad \forall \, n \in \mathcal{I}_M,
    \label{eq:wk}
\end{equation}
with $Y^{\sigma, \gamma}_m$ defined in Eq. \eqref{eq:Xs}. Direct summation of \eqref{eq:wk} costs $4Mh_{\gamma}$ flops, which we intend to improve. We now split the reference interval $I = [-1, 1]$ for a submatrix on level $\gamma$ into 2 smaller subintervals 
\begin{equation}
I(q) = \left[q-1, q\right], \quad q \in \{0, 1\},  \label{eq:subintervals}  
\end{equation}
such that $I = I(0)\cup I(1)$. The two subintervals are shown in Fig. \ref{fig:2levelsFMM}.
Next, we split the sum in Eq. \eqref{eq:wk} into $2$ parts corresponding to the 2 new subintervals
\begin{equation}
    {w}_n(\gamma, b, q) = \sum_{l\in \{0,1\}}\sum_{m=l h_{\gamma+1}}^{(l+1)h_{\gamma+1}-1} T_n (Y^{\sigma,\gamma}_m) {f}_m^{\sigma}(\gamma, b, q),
    \label{eq:wk2}
\end{equation}
where $Y^{\sigma,\gamma}_m \in I(l)$, for all $lh_{\gamma+1} \le m < (l+1)h_{\gamma+1}$. The objective now is to replace terms on the right hand side of Eq. \eqref{eq:wk2} with terms already computed on level $\gamma+1$. 
To this end we first recognise from Fig. \ref{fig:2levelsFMM} (or Eq. \eqref{eq:flevels}) that 
\begin{equation}
\{{f}_m^{\sigma}(\gamma+1, 2b+q+1, l)\}_{m=0}^{h_{\gamma+1}-1} = \{f_m^{\sigma}(\gamma, b, q)\}_{m=lh_{\gamma+1}}^{(l+1)h_{\gamma+1}-1}, \label{eq:f0}
\end{equation}
which can be used directly in Eq. \eqref{eq:wk2}. We get
\begin{equation}
    {w}_n(\gamma, b, q) = \sum_{l\in \{0,1\}}\sum_{m=0}^{h_{\gamma+1}-1} T_n (Y^{\sigma,\gamma}_{m+lh_{\gamma+1}}) {f}_m^{\sigma}(\gamma+1, 2b+q+1, l),
    \label{eq:wk22}
\end{equation}
where we see that for given $l$ the Chebyshev basis functions are using coordinates only for half the regular interval since $Y^{\sigma,\gamma}_{m+lh_{\gamma+1}} \in I(l)$ for $l \in \{0, 1\}$ and $m \in \mathcal{I}_{h_{\gamma+1}}$.

In order to map the half-interval Chebyshev basis functions in \eqref{eq:wk22} to the full-interval basis functions used on level $\gamma+1$ we introduce 
\begin{equation}
    T_n((Y+2q-1)/2) = \sum_{k=0}^{n} b^{(q)}_{nk} T_k(Y), \quad Y \in I, \quad q \in \{0, 1\}, \label{eq:Tky}
\end{equation}
where $(Y+2q-1)/2 \in I(q)$ for $q\in\{0, 1\}$. We can compute the matrices $\boldsymbol{B}(q) = (b^{q}_{nk})_{n,k=0}^{M-1}$ using the orthogonality of Chebyshev polynomials or, alternatively, from the binomial theorem and basis changes, as shown in Appendix \ref{sec:A3}. The matrices are lower triangular since the half-interval Chebyshev polynomials $T_n((Y+2q-1)/2)$ are still polynomials in $\mathbb{P}_n(I(q))$ and for the first 4 rows we get
\begin{equation}
    \boldsymbol{B}{(0)} = \begin{bmatrix}
        1 &  0 & 0 & 0 & 0 &\cdots \\
        -\frac{1}{2} & \frac{1}{2} & 0 & 0 & 0 & \cdots \\
        -\frac{1}{4} & -1 & \frac{1}{4} & 0 & 0 & \cdots \\
        \frac{1}{4} & \frac{3}{8} & -\frac{3}{4} & \frac{1}{8} & 0 & \cdots \\
        \vdots & \vdots &\vdots &\vdots &\vdots &\ddots 
    \end{bmatrix}, \quad 
    \boldsymbol{B}{(1)} = \begin{bmatrix}
        1 &  0 & 0 & 0 & 0 &\cdots \\
        \frac{1}{2} & \frac{1}{2} & 0 & 0 & 0 & \cdots \\
        -\frac{1}{4} & 1 & \frac{1}{4} & 0 & 0 & \cdots \\
        -\frac{1}{4} & \frac{3}{8} & \frac{3}{4} & \frac{1}{8} & 0 & \cdots \\
        \vdots & \vdots &\vdots &\vdots &\vdots &\ddots
    \end{bmatrix}.
    \label{eq:appT}
\end{equation}

Going back to Eq. \eqref{eq:wk22} we can now use Eq. \eqref{eq:Tky} on the right hand side, and since $Y^{\sigma, \gamma}_{m+lh_{\gamma+1}} = (Y^{\sigma, \gamma+1}_m+2l-1)/2$ for $m\in \mathcal{I}_{h_{\gamma+1}}$ and $l \in \{0, 1\}$ we have
\begin{equation}
    T_n(Y^{\sigma, \gamma}_{m+lh_{\gamma+1}}) = \sum_{k=0}^n b^{(l)}_{nk} T_k(Y^{\sigma, \gamma+1}_{m}), \quad \text{for } \, l \in \{0, 1\}, \, m \in \mathcal{I}_{h_{\gamma+1}}. \label{eq:Ttrunc}
\end{equation}
We obtain
\begin{equation}
    {w}_n(\gamma, b, q)  
    = \sum_{l \in (0, 1)} \sum_{k=0}^n b_{nk}^{(l)} \sum_{m=0}^{h_{\gamma+1}-1} T_k(Y^{\sigma,\gamma+1}_m) {f}_m^{\sigma}(\gamma+1, 2b+q+1, l),
    \label{eq:wk4}
\end{equation}
and realise that the last sum in Eq. \eqref{eq:wk4} has already been computed on level $\gamma+1$ with the result stored in $\boldsymbol{w}(\gamma+1, 2b+q+1, l)$. As such we can rearrange and simplify into
\begin{equation}
    {w}_n(\gamma, b, q) = \sum_{l\in \{0, 1\}}\sum_{k=0}^n b_{nk}^{(l)} {w}_k(\gamma+1, 2b+q+1, l),
    \label{eq:wk5}
\end{equation}
which in matrix form becomes Eq. \eqref{eq:wtk}.

\vskip 0.1in

\renewcommand{\thealgorithm}{a.\arabic{algorithm}}%
\setcounter{algorithm}{0}

\subsection{Computing Eq. \eqref{eq:wtk} efficiently}
\label{sec:A2}

Algorithm \ref{alg:Tw} computes Eq. \eqref{eq:wtk} in $M^2+2M$ operations using the fact that $\boldsymbol{B}(0)$ and $\boldsymbol{B}(1)$ differ only in the sign of the odd diagonals, see App. \ref{sec:A3}.

\begin{algorithm}
\caption{A fast computation of Eq. \eqref{eq:wtk} costing $M^2+2M$ flops, attributed to one lower triangular matrix-vector product plus setting up the work array $\boldsymbol{z}$.}
\label{alg:Tw}
\begin{algorithmic}[1]
\Function{eq43}{$\boldsymbol{w}^0, \boldsymbol{w}^1, \boldsymbol{B}$}
\Input{$\boldsymbol{w}^0$}{array $\in \mathbb{R}^{M}$}{$\boldsymbol{w}(\gamma,b,q)$ in Eq. \eqref{eq:wtk}}
\Input{$\boldsymbol{w}^1$}{array $\in \mathbb{R}^{2 \times M}$}{$\boldsymbol{w}(\gamma+1,2b+q+1)$ in Eq. \eqref{eq:wtk}}
\Input{$\boldsymbol{B}$}{$\mathbb{R}^{M \times M}$}{$\boldsymbol{B}(0)$ in Eq. \eqref{eq:wtk}}
\WorkArray{$\boldsymbol{z}$}{$\mathbb{R}^{2\times M}$}
\State $M \gets \mathbf{{len}}~\boldsymbol{w}^0$
\For{$i \gets 0, M-1$}
    \State $z_{0,i} \gets w^1_{0,i} + w^1_{1,i}$ 
    \State $z_{1,i} \gets w^1_{0,i} - w^1_{1,i}$ 
\EndFor
\State $n \gets 0$
\While {$n < M$}
    \State $k \gets n \,  \mathbf{mod} \, 2$
    \For{$i \gets 0, M-n-1$}
        \State $w^0_{i+n} \gets w^0_{i+n} + b_{i+n, i} z_{k,i}$
    \EndFor
    \State $n \gets n+1$
\EndWhile
\State \Return $\boldsymbol{w}^0$
\EndFunction
\end{algorithmic}
\end{algorithm}


\vskip 0.1in

\subsection{Computing the matrices $\boldsymbol{B}(l)$ in Eq. \eqref{eq:wtk} exactly}
\label{sec:A3}

We compute the matrices using a simple basis change and the binomial theorem. We have that
\begin{equation}
    x^n = \sum_{j=0}^n w_{nj}T_j(x), \quad n \in \mathcal{I}_M, \label{eq:xn}
\end{equation}
or with $\boldsymbol{\underline{x}} = (x^j)_{j=0}^{M-1}$, $\boldsymbol{W} = (w_{nj})_{n,j=0}^{M-1}$ and $\boldsymbol{\underline{T}} = (T_j(x))_{j=0}^{M-1}$
\begin{equation}
    \boldsymbol{\underline{x}} = \boldsymbol{W} \boldsymbol{\underline{T}}. \label{eq:WT}
\end{equation}
The lower triangular coefficient matrix $\boldsymbol{W}$, where only the even diagonals are nonzero, is given as \cite{mathar06}
\begin{equation}
w_{nj} =
    \begin{cases}
        \frac{2}{c_j2^n}\dbinom{n}{\frac{n-j}{2}}, \quad &n-j \text{ even and } j \le n, \\
        0, &\text{otherwise}.
    \end{cases}
\end{equation}
The inverse matrix $\boldsymbol{W}^{-1}$ is also lower triangular and given as
\begin{equation}
    w^{-1}_{nj} = \begin{cases}
        1, \quad &n=j=0, \\
        \frac{n(-1)^{\frac{n-j}{2}}2^j}{n+j}\dbinom{\frac{n+j}{2}}{\frac{n-j}{2}}, \quad &n-j \text{ even}, n > 0 \text{ and }  j \le n ,\\
        0, &\text{otherwise},
    \end{cases}
\end{equation}
rewritten here from \cite{mathar06} such that $\boldsymbol{\underline{T}} = \boldsymbol{W}^{-1}\boldsymbol{\underline{x}}$.

The matrices $\boldsymbol{B}(l)$ are used to express the half-interval Chebyshev polynomials $T_k((x+2l-1)/2)$ in terms of $\{T_n(x)\}_{n=0}^k$, where $x \in [-1,1]$ and $(x+2l-1)/2 \in [l-1, l]$ for $l\in\{0, 1\}$. The binomial theorem tells us that
\begin{equation}
    \left(\frac{x+2l-1}{2}\right)^n = \sum_{j=0}^n v^{(l)}_{nj}  x^j,
\end{equation}
where 
\begin{equation}
    v^{(l)}_{nj} = \begin{cases}
        \dbinom{n}{n-j}\frac{(2l-1)^{n-j}}{2^n}, \quad &j \le n, \\
        0, &\text{otherwise}.
    \end{cases} \label{eq:al}
\end{equation}
In matrix form, with $\boldsymbol{\overline{x}}{(l)} = \{((x+2l-1)/2)^j\}_{j=0}^{M-1}$ and $\boldsymbol{V}{(l)} = (v^{(l)}_{nj})_{n,j=0}^{M-1}$, we get
\begin{equation}
    \boldsymbol{\overline{x}} = \boldsymbol{V} \boldsymbol{\underline{x}}. \label{eq:AX}
\end{equation}
From Eq. \eqref{eq:xn} we have
\begin{equation}
    \boldsymbol{\overline{x}} = \boldsymbol{W} \overline{\boldsymbol{T}}, \label{eq:WTO} 
\end{equation}
where $ \overline{\boldsymbol{T}}(l) = \{T_j((x+2l-1)/2)\}_{j=0}^{M-1}$. Setting now Eq. \eqref{eq:AX} equal to \eqref{eq:WTO} and using Eq. \eqref{eq:WT} we get
\begin{equation}
    \boldsymbol{V} \boldsymbol{W} \boldsymbol{\underline{T}} = \boldsymbol{W} \boldsymbol{\overline{T}},
\end{equation}
and thus 
\begin{equation}
    \boldsymbol{\overline{T}} = \boldsymbol{W}^{-1}\boldsymbol{V} \boldsymbol{W} \boldsymbol{\underline{T}}. \label{eq:WAW}
\end{equation}
Comparing Eq. \eqref{eq:WAW} to Eq. \eqref{eq:Tky} we get that $\boldsymbol{B}(l) = \boldsymbol{W}^{-1}\boldsymbol{V}{(l)}\boldsymbol{W}$. From Eq. \eqref{eq:al}  it is obvious that $\boldsymbol{V}{(0)}$ and $\boldsymbol{V}{(1)}$ differ only in the sign of the odd diagonals. Since $\boldsymbol{W}$ and $\boldsymbol{W}^{-1}$ decouple the odd and even components (every odd diagonal is zero) this means that also $\boldsymbol{B}(0)$ and $\boldsymbol{B}(1)$ will differ only in the sign of the odd diagonals. Yet another interesting computational feature of the $\boldsymbol{B}(l)$ matrices is that for $M=18$ they can be represented exactly in double precision since all items are rational numbers with the denominator being a small power of 2.

\end{appendices}

\end{document}